\documentclass[11pt]{amsart}
\usepackage{amssymb, enumerate, xspace, graphicx}
\usepackage[all]{xy}
\SelectTips{cm}{}

\numberwithin{equation}{section}

1

\numberwithin{subsection}{section}

\allowdisplaybreaks[1]

\newenvironment{enumeratea}
{\begin{enumerate}[\upshape (a)]}
{\end{enumerate}}

\newtheorem*{namedtheorem}{\theoremname}
\newcommand{\theoremname}{testing}

\newtheorem{theorem}{Theorem}[section]
\newtheorem{proposition}[theorem]{Proposition}
\newtheorem{proposition-definition}[theorem]
{Proposition-Definition}
\newtheorem{corollary}[theorem]{Corollary}
\newtheorem{lemma}[theorem]{Lemma}

\theoremstyle{definition}

\newtheorem{remark}[theorem]{Remark}

\theoremstyle{remark}

\newcommand\nome{testing}
\newcommand\call[1]{\label{#1}\renewcommand\nome{#1}}
\newcommand\itemref[1]{\item\label{\nome;#1}}
\newcommand\refall[2]{\ref{#1}~(\ref{#1;#2})}
\newcommand\refpart[2]{(\ref{#1;#2})}

\newcommand\cA{\mathcal{A}}

\newcommand\cC{\mathcal{C}}

\newcommand\cO{\mathcal{O}}

\renewcommand\AA{\mathbb{A}}

\newcommand\CC{\mathbb{C}}

\newcommand\FF{\mathbb{F}}

\newcommand\PP{\mathbb{P}}
\newcommand\QQ{\mathbb{Q}}

\newcommand\ZZ{\mathbb{Z}}

\newcommand\rB{\mathrm{B}}
\newcommand\rC{\mathrm{C}}

\newcommand\rH{\mathrm{H}}
\newcommand\rI{\mathrm{I}}

\newcommand\rN{\mathrm{N}}
\newcommand\rO{\mathrm{O}}

\newcommand\rT{\mathrm{T}}

\newcommand\rme{\mathrm{e}}

\newcommand\rmi{\mathrm{i}}

\newcommand\arr{\ifinner\to\else\longrightarrow\fi}

\newcommand{\xarr}{\xrightarrow}

\newcommand\noqed{\renewcommand\qed{}}

\renewcommand\H{\operatorname{H}}

\newcommand\eqdef{\overset{\mathrm{\scriptscriptstyle def}} =}

\newcommand\into{\hookrightarrow}

\renewcommand\th{^\text{th}}

\def\displaytimes_#1{\mathrel{\mathop{\times}\limits_{#1}}}

\def\displayotimes_#1{\mathrel{\mathop{\bigotimes}\limits_{#1}}}

\renewcommand\hom{\operatorname{Hom}}

\newcommand\tor{\operatorname{Tor}}

\newcommand\spec{\operatorname{Spec}}

\newcommand\generate[1]{\langle #1 \rangle}

\newcommand\id{\mathrm{id}}

\newcommand\pr{\operatorname{pr}}

\newcommand\dash{\nobreakdash-\hspace{0pt}}

\newdir{ >}{{}*!/-5pt/@{>}}

\newcommand\doublelong[2]{\mathbin{\xymatrix{{}\ar@<3pt>[r]^{#1}
\ar@<-3pt>[r]_{#2}&}}}

\newlength{\ignora}
\newcommand{\hsmash}[1]{\settowidth{\ignora}{#1}#1\hspace{-\ignora}}

\newcommand\ch[1][*]{\operatorname{A}^{#1}_}

\renewcommand\H[1][*]{\operatorname{H}^{#1}_}

\renewcommand\c{\operatorname{c}}

\newcommand\sym{\operatorname{Sym}}

\newcommand\s[1][p]{{\mathrm{S}_{#1}}}

\newcommand\cyc[1][p]{{\mathrm{C}_{#1}}}

\newcommand\pgl[1][p]{{\mathrm{PGL}_{#1}}}

\newcommand\gm{\mathbb{G}_{\mathrm{m}}}

\newcommand\SL[1][p]{{\mathrm{SL}_{#1}}}

\newcommand\GL[1][p]{{\mathrm{GL}_{#1}}}

\newcommand\Sp[1][n]{{\mathrm{Sp}_{#1}}}

\newcommand\SO[1][n]{{\mathrm{SO}_{#1}}}

\newcommand\N[1][p]{{\mathrm{N}_{#1}}}

\newcommand\mmu[1][p]{{\boldsymbol \mu}_{#1}}

\newcommand\toruspgl[1][p]{{\mathrm{T}_{\pgl[#1]}}}

\newcommand\torusgl[1][p]{{\mathrm{T}_{\GL[#1]}}}

\newcommand\torussl[1][p]{{\mathrm{T}_{\SL[#1]}}}

\newcommand\chartorus[1]{\widehat{\mathrm{T}}_{#1}}

\newcommand\cycpgl[1][p]{{\cyc\ltimes\toruspgl[#1]}}

\newcommand\cycsl[1][p]{{\cyc\ltimes\torussl[#1]}}

\newcommand\cycgl[1][p]{{\cyc\ltimes\torusgl[#1]}}

\newcommand\sympgl[1][p]{{\s[#1]\ltimes\toruspgl[#1]}}

\newcommand\symgl[1][p]{{\s[#1]\ltimes\torusgl[#1]}}

\newcommand\norm[1][p]{\mathbb{F}_{#1}^*\ltimes\mathbb{F}_{#1}}

\newcommand\liesl[1][p]{\mathfrak{sl}_{#1}}

\newcommand\liegl[1][p]{\mathfrak{gl}_{#1}}

\newcommand\liesls[1][p]{\mathfrak{sl}_{#1}^0}

\newcommand\diag[1][p]{D_{#1}}

\newcommand\diags[1][p]{D^0_{#1}}

\newcommand\invert{\otimes\mathbb{Z}[1/(p-1)!]}

\newcommand\invertp{\otimes\mathbb{Z}[1/p]}

\newcommand\cycmu[1][p]{{\cyc[#1]\times \mmu[#1]}}

\newcommand\invtorus[1][p]{\bigl(\ch{\toruspgl[#1]}\bigr)^{\s[#1]}}

\newcommand{\slfinite}{\SL[2](\FF_{p})}

\newcommand\tors{_\mathrm{tors}}

\newcommand\chern{\operatorname{c}}

\newcommand\res[2]{\operatorname{res}^{#1}_{#2}}

\newcommand\tsf[2]{\operatorname{tsf}^{#1}_{#2}}

\newcommand\class{\mathrm{B}\,}

\newcommand{\even}{\mathrm{even}}
\newcommand{\odd}{\mathrm{odd}}

\begin{document}

\title[On the classifying space of $\pgl$]{On the cohomology and the Chow ring\\of the classifying space of $\pgl$}
\author{Angelo Vistoli}

\address{Dipartimento di Matematica\\
Universit\`a di Bologna\\
40126 Bologna\\ Italy}
\email{vistoli@dm.unibo.it}

\begin{abstract}
We investigate the integral cohomology ring and the Chow ring of the classifying space of the complex projective linear group $\pgl$, when $p$ is an odd prime. In particular, we determine their additive structures completely, and we reduce the problem of determining their multiplicative structures to a problem in invariant theory.
\end{abstract}

\date{April 24, 2005}

\thanks{Partially supported by the University of Bologna, funds for selected research topics, and by funds FIRB and PRIN from M.I.U.R..}

\maketitle

\tableofcontents


\section{Introduction}

Let $G$ be a complex linear group. One of the main invariants associated with $G$ is the cohomology $\H{G}$ of the classifying space $\mathrm{B}G$. B. Totaro (see \cite{totaro-classifying}) has also introduced an algebraic version of the cohomology of the classifying space of an algebraic group $G$ over a field $k$, the Chow ring $\ch{G}$ of the classifying space of $G$. When $k =  \CC$ there is a cycle ring homomorphism $\ch{G} \arr \H{G}$. Chow rings are normally infinitely harder to study than cohomology; it remarkable that, in contrast, $\ch{G}$ seems to be better behaved, and easier to study, than $\H{G}$. For example, when $G$ is a finite abelian group, $\ch{G}$ is the symmetric algebra over $\ZZ$ of the dual group $\widehat{G}$; while the cohomology ring contains this symmetric algebra, but is much more complicated (for example, will contain elements of odd degree), unless $G$ is cyclic.

This ring $\ch{G}$ has also been computed for $G = \GL[n]$, $\SL[n]$, $\Sp$ by Totaro (\cite{totaro-classifying}), for $G = \rO_{n}$ and $G = \SO[2n+1]$ by Totaro and R{.} Pandharipande (\cite{totaro-classifying} and \cite{rahul-orthogonal}), for $G = \SO[2n]$ by R{.} Field (\cite{field-so2n}), and for the semisimple simply connected group of type $G_{2}$ by N{.} Yagita (\cite{yagita-AHss}). Also, a lot of work has been done on the case of finite groups, by Totaro himself and by Guillot (\cite{guillot-chevalley} and \cite{guillot-steenrod}). However, not much was known for the $\pgl[n]$ series. Even the cohomology of $\rB\,\pgl[n]$ was mysterious. Algebraic topologists tend to work with cohomology with coefficients in a field, the case in which their extremely impressive toolkits work the best. When $p$ does not divide $n$, the cohomology ring $\H{}(\rB\,\pgl[n], \ZZ/p\ZZ)$ is a well understood polynomial ring. Also, since $\pgl[2] = \SO[3]$, the ring $\H{}(\rB\,\pgl[2], \ZZ/2\ZZ)$ is also well understood. The other results that I am aware of on $\H{}(\rB\,\pgl[n], \ZZ/p\ZZ)$ that were known before this article was posted are the following.

\begin{enumerate}

\item In \cite{kono-mimura-shimada}, the authors compute $\H{}(\rB\,\pgl[3], \ZZ/3\ZZ)$ as a ring, by presentations and relations.

\item The ring $\H{}(\rB\,\pgl[n], \ZZ/2\ZZ)$ is known when $n \equiv 2 \pmod 4$ (\cite{kono-mimura} and \cite{toda-classifying}).

\item Some interesting facts on $\H{}(\rB\,\pgl[p], \ZZ/p\ZZ)$ are proved in \cite{vavpetic-viruel}.

\end{enumerate}

On the other hand, to my knowledge no one had studied the integral cohomology ring $\H{\pgl[n]}$.

In the algebraic case, the only known results about $\ch{\pgl[n]}$, apart from the case of $\pgl[2] = \SO[3]$, concern $\pgl[3]$ and were proved by Vezzosi in \cite{vezzosi-pgl3}. Here he determines almost completely the structure of $\ch{\pgl[3]}$ by generators and relations; the only ambiguity is about one of the generators, denoted by $\chi$ and living in $\ch[6]{\pgl[3]}$, about which he knows that it is $3$-torsion, but is not able to determine whether it is $0$. This $\chi$ maps to $0$ in the cohomology ring $\H{\pgl[3]}$; according to a conjecture of Totaro, the cycle map $\ch{\pgl[3]} \arr \H{\pgl[3]}$ should be injective; so, if the conjecture is correct, $\chi$ should be $0$.

Despite this only partial success, the ideas in \cite{vezzosi-pgl3} are very important. The main one is to make use of the \emph{stratification method} to get generators. This is how it works. Recall that Edidin and Graham (\cite{edidin-graham-equivariant}) have generalized Totaro's ideas to give a full-fledged equivariant intersection theory. Let $V$ be a representation of a group $G$; then we have $\ch{G} = \ch{G}(V)$. Suppose that we have a stratification $V_{0}$, \dots,~$V_{t}$ of $V$ by locally closed invariant subvarieties, such that each $V_{\leq i} \eqdef \cap_{j\leq i}V_{j}$ is open in $V$, each $V_{i}$ is closed in $V_{\leq i}$, and $V_{t} = V \setminus\{0\}$. If we can determine generators for $\ch{G}(V_{i})$ for each $i$, then one can use the localization sequence
   \[
   \ch{G}(V_{i}) \arr \ch{G}(V_{\leq i}) \arr \ch{G}(V_{\leq i-1}) \arr 0
   \]
and induction to get generators for $\ch{G}(V \setminus\{0\})$; and since $\ch{G}(V \setminus\{0\}) = \ch{G}/\bigl(\chern_{r}(V)\bigr)$, where $\chern_{r}(V) \in \ch[r]{G}$ is the $r\th$ Chern class of $V$, we obtain that $\ch{G}$ is generated by lifts to $\ch{G}$ of the generators for $\ch{G}(V \setminus\{0\})$, plus $\chern_{r}(V)$.

The stratification methods gives a unified approach for all the known calculations of $\ch{G}$ for classical groups (see \cite{molina-vistoli-classifying}).

Vezzosi applies the method to the adjoint representation space $V = \liesl[3]$ consisting of matrices with trace $0$. The open subscheme $V_{0}$ is the subscheme of matrices with distinct eigenvalues; its Chow ring is related with the Chow ring of the normalizer $\rN_{3}$ of a maximal torus $\toruspgl$ in $\pgl[3]$. In order to get relations, Vezzosi uses an unpublished result of Totaro, implying that the restriction homomorphism $\ch{\pgl[3]} \arr \ch{\rN_{3}}$ is injective. The reason why he is not able to determine whether $\chi$ is $0$ or not is that he does not have a good description of the $3$-torsion in $\ch{\rN_{3}}$.

In this paper we refine Vezzosi's approach, and extend it to the case of $\pgl[p]$, where $p$ is an odd prime.

Let $\toruspgl$ be the standard maximal torus in $\pgl$, consisting of classes of diagonal matrices, $\rN_{p}$ its normalizer, $\s = \rN_{p}/\toruspgl$ its Weyl group. Here are our main results (see Section~\ref{sec:mainresults} for details).

\begin{enumerate}

\item The natural homomorphism $\ch{\toruspgl} \arr \invtorus$ is surjective, and has a natural splitting $\invtorus \arr \ch{\toruspgl}$, which is a ring homomorphism.

\item\label{item:ch} The ring $\ch{\pgl}$ is generated as an algebra over $\invtorus$ by a single $p$-torsion element $\rho \in \ch[p+1]{\pgl}$; we also describe the relations.

\item\label{item:H} The ring $\H{\pgl}$ is generated as an algebra over $\invtorus$ by two elements: the image $\rho \in \H[2p+2]{\pgl}$ of the class above and the Brauer class $\beta \in \H[3]{\pgl}$; we also describe the relations.

\item Using (\ref{item:ch}) and (\ref{item:H}) above, we describe completely the additive structures of $\ch{\pgl}$ and $\H{\pgl}$.

\item For $p = 3$ we give a presentation of $\invtorus[3]$ by generators and relations (this is already in \cite{vezzosi-pgl3}); and this, together with (\ref{item:ch}) and (\ref{item:H}) above, gives presentations of $\ch{\pgl[3]}$ and $\H{\pgl[3]}$, completing the work of \cite{vezzosi-pgl3}.

\item The cycle homomorphism $\ch{\pgl} \arr \H[\even]{\pgl}$ into the even-dimensional cohomology is an isomorphism.

\end{enumerate}

The ring $\invtorus$ is complicated when $p > 3$; see the discussion in Section~\ref{sec:invtorus}.

The class $\rho$ in (\ref{item:ch}) seems interesting, and gives a new invariant for sheaves of Azumaya algebras of prime rank (Remark~\ref{rmk:azumaya}). In \cite{targa}, Elisa Targa shows that $\rho$ is not a polynomial in Chern classes of representations of $\pgl$.

Many of the ideas in this paper come from \cite{vezzosi-pgl3}. The main new contributions here are the contents of Sections $6$ and $7$ (the heart of these results are Proposition~\ref{prop:chcycgl-hcycgl}, and the proof of Lemma~\ref{lem:keyfact-cycgl}), which substantially improve our understanding of the cohomology and Chow ring of the classifying space of $\rN_{p}$, and Proposition~\ref{prop:localization}, which gives a way of showing that in the stratification method no new generators come from the strata corresponding to non-zero matrices with multiple eigenvalues, thus avoiding the painful case-by-case analysis that was necessary in \cite{vezzosi-pgl3}.

Recently I received a preprint of M{.}~Komeko and N{.}~Yagita (\cite{kameko-yagita-pup}) who also calculate the additive structure of $\H{\pgl}$, with completely different methods.

\subsection*{Acknowledgments} I would like to thank Nitin Nitsure and Alejandro Adem for pointing out references \cite{kono-mimura}, \cite{toda-classifying} and \cite{vavpetic-viruel} to me.

I am also in debt with Alberto Molina, who discovered a serious mistake in the proof of Theorem~\ref{thm:main-splitting} given in a preliminary version of the paper, and with Marta Morigi, who helped me fix it.

Finally, I would like to acknowledge the very interesting discussions I have had with Nobuaki Yagita and Andrzej Weber on the subject of this paper.

\section{Notations and conventions}

All algebraic groups and schemes will be of finite type over a fixed field $k$ of charateristic $0$. Furthermore, we will fix an odd prime $p$, and assume that $k$ contains a fixed $p\th$ root of $1$, denoted by $\omega$.  When $k = \CC$, we take $\omega = \rme^{2\pi\rmi/p}$.

The hypothesis that the characteristic be $0$ is only used in the proof Theorem~\ref{thm:totaro-injectivity}, which should however hold over an arbitrary field. If so, it would be enough to assume here that the characteristic of $k$ be different from $p$.

Our main tool is Edidin and Graham's equivariant intersection theory (see \cite{edidin-graham-equivariant}), which works over an arbitrary field; when we discuss cohomology, instead, we will always assume that $k = \CC$. All finite groups will be considered as algebraic groups over $k$, in the usual fashion. We denote by $\gm$ the multiplicative group of non-zero scalars over $k$, $\mmu[n]$ the algebraic group of $n\th$ roots of $1$ over $k$.

Whenever $V$ is a vector space over $k$, we also consider it as a scheme over $k$, as the spectrum of the symmetric algebra of the dual vector space $V^{\vee}$. If $V$ is a representation of an algebraic group $G$, then there is an action of $G$ on $V$ as a scheme over $k$.

We denote by $\torusgl$, $\torussl$ and $\toruspgl$ the standard maximal tori in the respective groups, those consisting of diagonal matrices. We identify the Weyl groups of these three groups with the symmetric group $\s$. We also denote the normalizer of $\toruspgl$ in $\pgl$ by $\sympgl$.

If $a_{1}$, \dots,~$a_{p}$ are elements of $k^{*}$, we will denote by $[a_{1}, \dots, a_{p}]$ the diagonal matrix in $\GL$ with entries $a_{1}$, \dots,~$a_{p}$, and also its class in $\pgl$. In general, we will ofter use the same symbol for a matrix in $\GL$ and its class in $\pgl$; this should not give rise to confusion.

It is well known that the arrows
   \[
   \ch{\GL} \longrightarrow (\ch{\torusgl})^{\s}, \quad
   \H{\GL} \longrightarrow (\ch{\torusgl})^{\s}
   \]
and
   \[
   \ch{\SL} \longrightarrow (\ch{\torussl})^{\s}, \quad
   \H{\GL} \longrightarrow (\ch{\torussl})^{\s}
   \]
induced by the embeddings $\torusgl \into \GL$ and $\torussl \into \SL$ are isomorphisms. If we denote by $x_i \in \ch{\torusgl} = \H{\torusgl}$ the first Chern class of the $i\th$ projection $\torusgl \to \gm$, or its restriction to $\torussl$, then $\ch{\torusgl} = \H{\torusgl}$ is the polynomial ring $\mathbb{Z}[x_1, \dots, x_p]$, while $\ch{\torussl} = \H{\torussl}$ equals  $\mathbb{Z}[x_1, \dots, x_p]/(x_1 + \dots  + x_p)$. If we denote by $\sigma_1$, \dots,~$\sigma_p$ the elementary symmetric functions in the $x_i$, then we conclude that
   \[
   \ch{\GL} = \H{\GL} = \mathbb{Z}[\sigma_1, \dots, \sigma_p]
   \]
while
   \[
   \ch{\SL} = \H{\SL} = \mathbb{Z}[\sigma_1, \dots, \sigma_p]/(\sigma_1)
   = \mathbb{Z}[\sigma_2, \dots, \sigma_p].
   \]

The ring $\ch{\toruspgl} = \H{\toruspgl}$ is the subring of $\ch{\torusgl}$ generated by the differences $x_i - x_j$. In particular it contains the element $\delta = \prod_{i \neq j} (x_i - x_j)$, which we call the \emph{discriminant} (up to sign, it is the classical discriminant); it will play an important role in what follows.

We will use the following notation: if $R$ is a ring, $t_{1}$, \dots,~$t_{n}$ are elements of $R$, $f_{1}$, \dots,~$f_{r}$ are polynomials in $\ZZ[x_{1}, \dots, x_{n}]$, we write \[ R = \ZZ[t_{1}, \dots, t_{n}]/ \bigl(f_{1}(t_{1}, \dots, t_{n}), \dots, f_{r}(t_{1}, \dots, t_{n})\bigr) \] to indicate the the ring $R$ is generated by $t_{1}$, \dots,~$t_{n}$, and the kernel of the evaluation map $\ZZ[x_{1}, \dots, x_{n}] \arr R$ sending $x_{i}$ to $t_{i}$ is generated by $f_{1}$, \dots,~$f_{r}$. When there are no $f_{i}$ this means that $R$ is a polynomial ring in the $t_{i}$.

\section{The main results}\label{sec:mainresults}

Consider the embedding $\mmu \into \toruspgl$ defined by $\zeta \mapsto [\zeta, \zeta^{2}, \dots, \zeta^{p-1}, 1]$. This induces a restriction homomorphism
   \[
   \ch{\toruspgl} \to \ch{\mmu} = \mathbb{Z}[\eta]/(p\eta),
   \]
where $\eta$ is the first Chern class of the embedding $\mmu \subseteq \gm$.

The restriction of the discriminant $\delta \in (\ch[p^2-p]{\toruspgl})^{\s}$ to $\mmu$ is the element
   \[
   \prod_{i \neq j} (i\eta - j\eta) =
   \biggl(\prod_{i \neq j} (i - j)\biggr) \eta^{p^2-p}
   \]
of $\mathbb{Z}[\eta]/(p\eta)$; this is non-zero multiple of $\eta^{p^2-p}$ (in fact, it is easy to check that it equals $-\eta^{p^2-p}$).

\begin{proposition}\label{prop:ker-whole}
The image of the restriction homomorphism
   \[
   \invtorus \longrightarrow \mathbb{Z}[\eta]/(p\eta)
   \]
is the subring generated by $\eta^{p^2-p}$.
\end{proposition}

This is proved at the end of Section~\ref{sec:cycpgl}.

\begin{theorem}\label{thm:main-splitting}
There exists a canonical ring homomorphism
   \[
   \invtorus \to \ch{\pgl}
   \]
whose composite with the restriction homomorphism $\ch{\pgl} \to \invtorus$ is the identity.
\end{theorem}

This is proved in Section~\ref{sec:splitting}.

As a consequence, $\ch{\pgl}$ and $\H{\pgl}$ can be regarded as $\invtorus$-algebras.

\begin{theorem}\label{thm:main-ch}
The $\invtorus$-algebra $\ch{\pgl}$ is generated by an element $\rho \in
\ch[p+1]{\pgl}$, and the ideal of relations is generated by the following:

\begin{enumeratea}

\item $p \rho = 0$, and

\item $\rho u = 0$ for all $u$ in the kernel of the homomorphism
$\invtorus \to \ch{\mmu}$.

\end{enumeratea}

\end{theorem}

There is a similar description for the cohomology: besides the element $\rho$, now considered as living in $\H[2p+2]{\pgl}$, we need a single class $\beta$ in degree $3$. This class is essentially the tautological \emph{Brauer class}. That is, if we call $\mathcal{C}$ the sheaf of complex valued continuous functions and $\mathcal{C}^*$ the sheaf of complex valued nowhere vanishing continuous functions on the classifying space $\class\pgl$, the tautological $\pgl$ principal bundle on $\class\pgl$ has a class  in the topological Brauer group $\H[2]{}(\class\pgl, \mathcal{C}^*)\tors$ (see \cite{grothendieck-brauer1}). On the other hand, the exponential sequence
   \[
   0 \longrightarrow \mathbb{Z} \xarr{2\pi\rmi}
   \mathcal{C} \longrightarrow \mathcal{C}^* \longrightarrow 1
   \]
induces a boundary homomorphism
   \[
   \H[2]{}(\class\pgl, \mathcal{C}^*) \longrightarrow
   \H[3]{}(\class\pgl, \mathbb{Z}) = \H[3]{\pgl},
   \]
which is an isomorphism, since $\class\pgl$ is paracompact, hence 
   \[
   \H[i]{}(\class\pgl, \mathcal{C}) = 0
   \]
for all $i > 0$. Our class $\beta$ is, up to sign, the image under this boundary homomorphism of the Brauer class of the tautological bundle. 

\begin{theorem}\label{thm:main-H}
The ring\/ $\H{\pgl}$ is the commutative $\invtorus$-algebra generated by an element $\beta$ of degree $3$ and the element $\rho$ of degree $2p+2$. The ideal of relations is generated by the following:

\begin{enumeratea}

\item $\beta^2 = 0$,

\item $p \rho = p \beta = 0$, and

\item $\rho u = \beta u = 0$ for all $u$ in the kernel of the homomorphism $\invtorus \to \ch{\mmu}$.

\end{enumeratea}

\end{theorem}

\begin{corollary}\label{cor:isom-even}
The cycle homomorphism induces an isomorphism of $\ch{\pgl}$ with $\H[\even]{\pgl}$.
\end{corollary}

From here it is not hard to get the additive structure of $\ch{\pgl}$ and $\H{\pgl}$. 
For each integer $m$, denote by $r(m,p)$ the number of partitions of $m$ into numbers between $2$ and $p$. If we denote by $\pi(m,p)$ the number of partitions of $m$ with numbers at most equal to $p$ (a more usual notation for this is $p(m,p)$, which does not look very good), then $r(m,p) = \pi(m,p) - \pi(m-1,p)$.

We will also denote by $s(m,p)$ the number of ways of writing $m$ as a linear combination $(p^{2} - p)i + (p+1)j$, with $i \geq 0$ and $j > 0$; and by $s'(m,p)$ the number of ways of writing $m$ as a the same linear combination, with $i \geq 0$ and $j \geq 0$. Obviously we have $s'(m,p) = s(m,p)$, unless $m$ is divisible by $p^{2}-p$, in which case $s'(m,p) = s(m,p) + 1$.

\begin{theorem} \label{thm:additive-structure}\hfil
\begin{enumeratea}

\item The groups $\ch[m]{\pgl}$ is isomorphic to
   \[
   \ZZ^{r(m,p)} \oplus (\ZZ/p\ZZ)^{s(m,p)}.
   \]

\item The group $\H[m]{\pgl}$ is isomorphic to $\ch[m/2]{\pgl}$ when $m$ is even, and is isomorphic to
   \[
   (\ZZ/p\ZZ)^{s'\left(\frac{m-3}{2},p\right)}
   \]
when $m$ is odd.
\end{enumeratea}
\end{theorem}

When $p =3$ we are able to get a description of $\ch{\pgl[3]}$ and $\H{\pgl[3]}$ by generator and relations, completing the work of \cite{vezzosi-pgl3}.

\begin{theorem}\label{thm:pgl3}\hfil

\begin{enumeratea}

\item $\ch{\pgl[3]}$ is the commutative $\ZZ$-algebra generated by elements $\gamma_{2}$, $\gamma_{3}$, $\delta$, $\rho$, of degrees $2$, $3$, $6$ and $4$ respectively, with relations
   \[
   27\delta - (4\gamma_{2}^{3} + \gamma_{3}^{2}), \quad
   3\rho, \quad
   \gamma_{2}\rho,\quad
   \gamma_{3}\rho.
   \]

\item $\H{\pgl[3]}$ is the commutative $\ZZ$-algebra generated by elements $\gamma_{2}$, $\gamma_{3}$, $\delta$, $\rho$ and $\beta$ of degrees $4$, $6$, $12$, $8$ and $3$ respectively, with relations
   \[
   27\delta - (4\gamma_{2}^{3} + \gamma_{3}^{2}), \quad
   3\rho, \quad 3\beta, \quad
   \beta^{2}, \quad
   \gamma_{2}\rho,\quad
   \gamma_{3}\rho, \quad
   \gamma_{2}\beta, \quad
   \gamma_{3}\beta.
   \]

\end{enumeratea}

\end{theorem}

The rest of the paper is dedicated to the proofs of these results. We start by recalling some basic facts on equivariant intersection theory.

\section{Preliminaries on equivariant intersection theory}\label{sec:preliminaries}

In this section the base field $k$ will be arbitrary.

We refer to \cite{totaro-classifying}, \cite{edidin-graham-equivariant} and \cite{vezzosi-pgl3} for the definitions and the basic properties of the Chow ring $\ch{G}$ of the classifying space of an algebraic group $G$ over a field $k$, and of the Chow group $\ch{G}(X)$ when $X$ is a scheme, or algebraic space, over $k$ on which $G$ acts, and their main properties. Almost all $X$ that appear in this paper will be smooth, in which case $\ch{G}(X)$ is a commutative ring; the single exception will be in the proof of Lemma~\ref{lem:keyfact-cycgl}. 

The connection between these two notions is that $\ch{G} = \ch{G}(\spec k)$.

Recall that $\ch{G}(X)$ is contravariant for equivariant morphism of smooth varieties; that is, if $f\colon X \arr Y$ is a $G$-equivariant morphism of smooth $G$-schemes, there in an induced ring homomorphism $f^{*}\colon \ch{G}(X) \arr \ch{G}(Y)$.

If $k = \CC$, and $X$ is a smooth algebraic variety on which $G$ acts, there is a cycle ring homomorphism $\ch{G}(X) \arr \H{G}(X)$ from the equivariant Chow ring to the equivariant cohomology ring; this is compatible with pullbacks. 

Furthermore, if $f$ is proper there is a pushforward $f_{*}\colon \ch{G}(Y) \arr \ch{G}(X)$; this is not a ring homomorphism, but it satisfies the projection formula
   \[
   f_{*}(\xi\cdot f^{*}\eta) = f_{*}\xi \cdot \eta
   \]
for any $\xi \in \ch{G}(X)$ and $\eta \in \ch{G}(Y)$.

When $Y$ is a closed $G$-invariant subscheme of $X$ and we denote by $\iota\colon Y \into X$ the embedding, then we have a localization sequence
   \[
   \ch{G}(Y) \stackrel{\iota_{*}}{\longrightarrow}
   \ch{G}(X) \arr \ch{G}(X \setminus Y) \arr 0.
   \]
The analogous statement for cohomology is different: here the restriction homomorphism $\H{G}(X) \arr \H{G}(X \setminus Y)$ is not necessarily surjective. However, when $X$ and $Y$ are smooth we have a long exact sequence
   \[
   \xymatrix@C+20pt{
   &\cdots\ar[r]& \H[i-1]{G}(X \setminus Y) \ar[dll] _{\partial} \\
   \H[i-2r]{G}(Y) \ar[r]_{\iota_{*}} 
   & \H[i]{G}(X) \ar[r]
   &\H[i]{G}(X \setminus Y) \ar[dll]_{\partial}\\
   \H[i -2r + 1]{G}(Y) \ar[r] & \cdots
   }
   \]
where $r$ is the codimension of $Y$ in $X$.

Furthermore, if $H \arr G$ is a homomorphism of algebraic groups, and $G$ acts on a smooth scheme $X$, we can define an action of $H$ on $X$ by composing with the given homomorphism $H \arr G$. Then we have a restriction homomorphism
   \[
   \res{G}{H}\colon \ch{G}(X)\arr \ch{H}(X).
   \]

Here is another property that will be used often. Suppose that $H$ is an algebraic subgroup of $G$. We can define a ring homomorphism $\ch{G}(G/H) \arr \ch{H}$ by composing the restriction homomorphism $\ch{G}(G/H) \arr \ch{H}(G/H)$ with the pullback $\ch{H}(G/H) \arr \ch{H}(\spec k) = \ch{H}$ obtained by the homomorphism $\spec k \arr G/H$ whose image is the image of the identity in $G(k)$. Then this ring homomorphism is an isomorphism.

More generally, suppose that $H$ acts on a scheme $X$. We define the induced space $G \times^{H} X$ as usual, as the quotient $(G \times X)/H$ by the free right action given by the formula $(g,x)h = (gh, h^{-1}x)$. This carries a natural left action of $G$ defined by the formula $g'(g,x) = (g'g,x)$. There is also an embedding $X \simeq H \times^{H}X \into G \times^{H} X$ that is $H$-equivariant: and the composite of the restriction homomorphism $\ch{G}(G \times^{H} X) \arr \ch{H}(G \times^{H} X)$ with the pullback $\ch{H}(G \times^{H} X) \arr \ch{H}(X)$ is an isomorphism.

Furthermore, if $V$ is a representation of $G$, then there are Chern classes $\chern_{i}(V)\in \ch[i]{G}$, satisfying the usual properties. More generally, if $X$ is a smooth scheme over $k$ with an action of $G$, every $G$-equivariant vector bundles $E \arr X$ has Chern classes $\chern_{i}(E) \in \ch[i]{G}(X)$.

The following fact will be used often.

\begin{lemma}\label{lem:exact-top}
Let $E \arr X$ be an equivariant vector bundle of constant rank $r$, $s\colon X \arr E$ the $0$-section, $E_{0} \subseteq E$ the complement of the $0$-section. Then the sequence
   \[
   \ch{G}(X) \xarr{\chern_{r}(E)} \ch{G}(X) \arr \ch{G}(E_{0}) \arr 0,
   \]
where the second arrow is the pullback along $E_{0} \arr \spec k$, is exact.

Furthermore, when $k = \CC$ we also have a long exact sequence
   \[
   \xymatrix@C+20pt{
   &\cdots\ar[r]& \H[i-1]{G}(E_{0}) \ar[dll] _{\partial} \\
   \H[i-2r]{G}(X) \ar[r]_{\chern_{r}(E)} 
   & \H[i]{G}(X) \ar[r]
   &\H[i]{G}(E_{0}) \ar[dll]_{\partial}\\
   \H[i -2r + 1]{G}(X) \ar[r] & \cdots
   }
   \]

\end{lemma}

\begin{proof}
This follows from the following facts:

\begin{enumerate}

\item the pullbacks $\ch{G}(X) \arr \ch{G}(E)$ and $\H{G}(X) \arr \H{G}(E)$ are isomorphisms,

\item the self-intersection formula, that says that the homomorphisms $s^{*}s_{*}\colon \ch{G}(X)\arr \ch{G}(X)$ and $s^{*}s_{*}\colon \H{G}(X)\arr \H{G}(X)$ are multiplication by $\chern_{r}(E)$, and

\item the localization sequences for Chow rings and cohomology. \qedhere

\end{enumerate}
\end{proof}

Let us recall the following results from \cite{totaro-classifying}.

\begin{enumerate}

\item If $T = \gm^{n}$ is a torus, and we denote by $x_{i} \in \ch[1]{T}$ the first Chern class of the  $i\th$ projection $T \arr \gm$, considered as a representation, then
   \[
   \ch{T} = \ZZ[x_{1}, \dots, x_{n}].
   \]
   
\item If $\torusgl[n]$ is the standard maximal torus in $\GL[n]$ consisting of diagonal matrices, then the restriction homomorphism $\ch{\GL[n]} \arr \ch{\torusgl[n]}$ induces an isomorphism
   \begin{align*}
   \ch{\GL[n]} &\simeq \ZZ[x_{1}, \dots, x_{n}]^{\s[n]}\\
   &= \ZZ[\sigma_{1}, \dots, \sigma_{n}]
   \end{align*}where the $\sigma_{i}$ are the elementary symmetric functions of the $x_{i}$.

\item If $\torussl[n]$ is the standard maximal torus in $\SL[n]$ consisting of diagonal matrices, and we denote by $x_{i}$ the restriction to $\ch{\SL[n]}$ of $x_{i} \in \ch{\GL[n]}$, then we have
   \[
   \ch{\torussl[n]} = \ZZ[x_{1}, \dots, x_{n}]/(\sigma_{1});
   \]
furthermore the restriction homomorphism $\ch{\SL[n]} \arr \ch{\torussl[n]}$ induces an isomorphism
   \begin{align*}
   \ch{\SL[n]} &\simeq
   \bigl(\ZZ[x_{1}, \dots, x_{n}]/(\sigma_{1})\bigr)^{\s[n]}\\
   &= \ZZ[\sigma_{1}, \sigma_{2}, \dots, \sigma_{n}]/(\sigma_{1}).
   \end{align*}
   
\item If $t \in \ch{\mmu[n]}$ is the first Chern class of the embedding $\mmu[n] \into \gm$, considered as a 1-dimensional representation, then we have
   \[
   \ch{\mmu[n]} = \ZZ[t](nt).
   \]
\end{enumerate}

Furthemore, if $G$ is any of the groups above and $k = \CC$, then the cycle homomorphism $\ch{G} \arr \H{G}$ is an isomorphism.

The following result is implicit in \cite{totaro-classifying}. Let $G$ be a finite algebraic group that is a product of copies of $\mmu[n]$, for various $n$. This is equivalent to saying that $G$ is a finite diagonalizable group scheme, or that $G$ is the Cartier dual of a finite abelian group $\Gamma$, considered as a group scheme over $k$. By Cartier duality, we have that $\Gamma$ is the character group $\widehat{G} \eqdef \hom(G, \gm)$.

\begin{proposition}
Consider the group homomorphism $\widehat{G} \arr \ch[1]{G}$ that sends each character $\chi\colon G \arr \gm$ into $\chern_{1}(\chi)$. The induced ring homomorphism $\sym_{\ZZ}\widehat{G} \arr \ch{G}$ is an isomorphism.
\end{proposition}

A more concrete way of stating this is the following. Set
   \[
   G = \mmu[n_{1}] \times \dots \times \mmu[n_{r}].
   \]
For each $i = 1$, \dots,~$n$ call $\chi_{i}$ the character obtained by composing the $i\th$ projection $G \arr \mmu[n_{i}]$ with the embedding $\mmu[n_{i}] \into \gm$, and set $\xi_{i} = \chern_{1}(\chi_{i}) \in \ch[1]{G}$. Then
   \[
   \ch{G} = \ZZ[\xi_{1}, \dots, \xi_{r}]/(n_{1}\xi_{1}, \dots, n_{r}\xi_{r}).
   \]

\begin{proof}
When $G = \mmu[n]$, this follows from Totaro's calculation cited above. The general case follows by induction on $r$ from the following Lemma.

\begin{lemma}
If $H$ is a linear algebraic group over $k$, the ring homomorphism
   \[
   \ch{H}\otimes_{\ZZ} \ch{\mmu[n]} \arr \ch{H \times\mmu[n]}
   \]
induced by the pullbacks $\ch{H} \arr \ch{H \times\mmu[n]}$ and $\ch{\mmu[n]} \arr \ch{H \times\mmu[n]}$ along the two projections $H \times\mmu[n] \arr H$ and $H \times\mmu[n] \arr \mmu[n]$ is an isomorphism.
\end{lemma}

\begin{proof}
This follows easily, for example, from the Chow--K\"unneth formula in \cite[Section~6]{totaro-classifying}, because $\mmu$ has a representation $V = k^{n}$ on which it acts by multiplication, with an open subscheme $U \eqdef V \setminus \{0\}$ on which it acts trivally; and the quotient $U/\mmu$ is the total space of a $\gm$-torsor on $\PP^{n-1}$, and, as such, it is a union of open subschemes of affine spaces.

It is also not hard to prove directly, as in \cite{molina-vistoli-classifying}.
\end{proof}
\noqed
\end{proof}

There is also a very important \emph{transfer} operation (sometimes called \emph{induction}). Suppose that $H$ is an algebraic subgroup of $G$ of finite index. The transfer homomorphism
   \[
   \tsf{H}{G}\colon  \ch{H} \arr  \ch{G}
   \]
(see \cite{vezzosi-pgl3}) is the proper pushforward from $\ch{H}\simeq \ch{G}(G/H)$ to $\ch{G}(\spec k) = \ch{G}$.

This is not a ring homomorphism; however, the projection formula holds, that is, if $\xi \in \ch{G}(X)$ and $\eta \in \ch{H}(X)$, we have
   \[
   \tsf{H}{G}\bigl(\xi \cdot  \res{G}{H}\eta\bigr) = \xi \cdot \tsf{H}{G}\eta
   \]
(in other words, $\tsf{H}{G}$ is a homomorphism of $\ch{G}(X)$-modules).

We are going to need the analogue of Mackey's formula in this context. Let $H$ and $K$ be algebraic subgroups of $G$, and assume that $H$ has finite index in $G$. We will also assume that the quotient $G/H$ is reduced, and a disjoint union of copies of $\spec k$ (this is automatically verified when $k$ is algebraically closed of characteristic~$0$). Then it is easy to see that the double quotient $K \backslash G /H$ is also the disjoint union of copies of $\spec k$. Furthermore, we assume that every element of $(K \backslash G /H)(k)$ is in the image of some element of $G(k)$. Of course this will always happen if $k$ is algebraically closed; with some work, this hypothesis can be removed, but it is going to be verified in all the cases to which we will apply the formula).

Denote by $\cC$ a set of representatives in $G(k)$ for classes in $(K \backslash G /H)(k)$. For each $s \in \cC$, set
   \[
   K_{s} \eqdef K \cap sHs^{-1} \subseteq G.
   \]
Obviously $K_{s}$ is a subgroup of finite index of $K$; there is also an embedding $K_{s} \into H$ defined by $k \mapsto s^{-1}ks$.

\begin{proposition}[Mackey's formula] \label{prop:mackey}\hfil
   \[
   \res{G}{K}\tsf{H}{G} = \sum_{s \in \cC} \tsf{K_{s}}{K}\res{H}{K_{s}}
   \colon \ch{H} \arr \ch{K}.
   \]
\end{proposition}

\begin{proof}
This is standard. We have that the equivariant cohomology rings $\ch{G}(G/H)$ and $\ch{G}(G/K)$ are canonically isomorphic to $\ch{H}$ and $\ch{K}$, respectively. The retriction homomorphism $\ch{G} \arr \ch{K}$ corresponds to the pullback $\ch{G}(\spec k) \arr \ch{G}(G/K)$, and the tranfer homomorphism corresponds to the proper pushforward $\ch{G}(G/H) \arr \ch{G}(\spec k)$. 

Since proper pushforwards and flat pullbacks commute, from the cartesian diagram
   \[
   \xymatrix{
   G/K \times G/H \ar[r]^-{\pr_{2}} \ar[d]^-{\pr_{1}}
   & G/H \ar[d]^{\pi}\\
   G/K \ar[r]^{\rho}
   & \spec k
   }
   \]
we get the equality
   \[
   \res{G}{K}\tsf{H}{G} = \rho^{*}\pi_{*} = 
   \pr_{1*}\pr_{2}^{*}\colon \ch{H} \arr \ch{K}.
   \]

Now we need to express $G/K \times G/H$ as a disjoint union of orbits by the diagonal action of $G$. There is a $G$-invariant morphism $G \times G \arr G$, defined by the rule $(a,b) \mapsto a^{-1}b$, that induces a morphism $G/K \times G/H \arr K \backslash G /H$. For each $s \in \cC$, call $\Omega_{s}$ the inverse image of $s \in (K \backslash G /H)(k)$, so that $G/K \times G/H$ is a disjoint union $\coprod_{s \in \cC} \Omega_{s}$. It is easy to verify that $\Omega_{s}$ is the orbit of the class $[1,s] \in (G/K \times G/H)(k)$ of the element $(1,s) \in (G \times G)(k)$, and that the stabilizer of $[1,s]$ is precisely $K_{s}$. From this we get an isomorphism
   \[
   G/K \times G/H \simeq \coprod_{s \in \cC} G/K_{s}
   \]
from which the statement follows easily.
\end{proof}

\begin{proposition}\label{prop:envelope}
Assume that $G$ is smooth. Let $f \colon X \arr Y$ a proper $G$-equivariant morphism of $G$-schemes. Assume that for every $G$-invariant closed subvariety $W \subseteq Y$ there exists a $G$-invariant closed subvariety of $X$ mapping birationally onto $W$. Then the pushforward $f_{*} \colon \ch{G}X \arr \ch{G}Y$ is surjective.
\end{proposition}

Here by $G$-invariant closed subvariety of $X$ we mean a closed subscheme $V$ of $X$ that is reduced, and such that $G$ permutes the irreducible components of $V$ transitively (one sometimes says that $V$ is \emph{primitive}).

This property can be expressed by saying that $X$ is an \emph{equivariant Chow envelope of $Y$} (see \cite[Definition~18.3]{fulton}).

\begin{proof}
In the non-equivariant setting the result follows from the definition of proper pushforward.

In our setting, let us notice first of all that if $Y' \arr Y$ is a $G$-equivariant morphism and $X' \eqdef Y' \times_{Y} X$, the projection $X' \arr Y'$ is also an equivariant Chow envelope (this is easy, and left to the reader). Therefore, if $U$ is an open subscheme of a representation of $G$ on which $G$ acts freely, the morphism $f\times\id_{U}\colon X \times U \arr Y \times U$ is an equivariant  Chow envelope. But since $G$ is smooth, it is easily seen that pullback from $(X \times U)/G$ to $X\times U$ defines a bijective correspondence between closed subvarieties of $(X \times U)/G$ and closed invariant subvarieties of $X\times U$; hence the $(X \times U)/G$ is a Chow envelope of $(Y \times U)/G$. So the proper pushforward $\ch{}\bigl((X \times U)/G\bigr) \arr \ch{}\bigl((Y \times U)/G\bigr)$ is surjective, and this completes the proof.
\end{proof}

\section{On $\cycmu$} \label{sec:cycmu}

A key role in our proof is played by a finite subgroups $\cycmu \subseteq \pgl$. 

We denote by $\cyc \subseteq \s$ the subgroup generated by the cycle $\sigma \eqdef (1\, 2\, \dots\,p)$. We embed $\s$ into $\pgl$ as usual by identifying a permutation $\alpha \in \s$ with the corresponding permutation matrix, obtained by applying $\alpha$ to the indices of the canonical basis $\rme_{1}$, \dots,~$\rme_{p}$ of $V$ (so that $\sigma \rme_{i} = \rme_{i+1}$, where addition is modulo $p$).

If we denote by $\tau$ the generator
   \[
   [\omega, \dots, \omega^{p-1}, 1]
   \]
of $\mmu \subseteq \pgl$,
we have that
   \[
   \tau\sigma = \omega\sigma\tau \quad \text{in }\GL;
   \]
so $\sigma$ and $\tau$ commute in $\pgl$, and they generate a subgroup
   \[
   \cycmu \subseteq \pgl.
   \]
We denote by $\alpha$ and $\beta$ the characters $\cycmu \arr \gm$ defined as
   \[
   \alpha(\sigma) = \omega \quad \text{and}\quad \alpha(\tau) = 1
   \]
and
   \[
   \beta(\sigma) = 1 \quad \text{and}\quad \beta(\tau) = \omega.
   \]

The following fact will be useful later.

\begin{lemma}\label{lem:describe-pgl}
If $i$ and $j$ are integers between $1$ and $p$, consider the matrix $\sigma^{i}\tau^{j}$ in the algebra $\liegl$ of $p \times p$ matrices. Then if $(i,j) \neq (p,p)$, the matrix $\sigma^{i}\tau^{j}$ has trace $0$, and its eigenvalues are precisely the $p$-roots of $1$.

Each $\sigma^{i}\tau^{j}$ is a semi-invariant for the action of $\cycmu$, with character $\alpha^{-j}\beta^{i}$. Furthermore the $\sigma^{i}\tau^{j}$ form a basis of $\liegl$, and those with $(i,j) \neq (p,p)$ form a bases of $\liesl$.

\end{lemma}

\begin{proof}
The fact that $\cycmu$ acts on $\sigma^{i}\tau^{j}$ via the character $\alpha^{-j}\beta^{i}$ is an elementary calculation, using the relation $\tau\sigma = \omega\sigma\tau$. From this it follows that the $\sigma^{i}\tau^{j}$ are linearly independent, and therefore form a basis of $\liegl$. The statement about the trace is also easy.

Let us check that the $\sigma^{i}\tau^{j}$ with $(i,j) \neq (p,p)$ have the elements of $\mmu$ as eigenvalues. When $i = p$ we get a diagonal matrix with eigenvalues are $\omega^{j}$, \dots,~$\omega^{pj}$, which are all the elements of $\mmu$, because $p$ is a prime and $j$ is not divisible by $p$. Assume that $i \neq p$. The numbers $i$, $2i$, \dots,~$pi$, reduced modulo $p$, coincide with $1$, \dots,~$p$. If $\lambda$ is a $p\th$ root of $1$, and $\rme_{1}$, \dots,~$\rme_{p}$ is the canonical basis of $k^{n}$, then the vector 
   \[
   \sum_{t=1}^{p}\lambda^{-t}\omega^{ij\binom{t}{2}}\rme_{ti}
   \]
is easily seen to be an eigenvector of $\sigma^{i}\tau^{j}$ with eigenvalue $\lambda$ (using the fact that
   \[
   \binom{t_{1}}{2} \equiv \binom{t_{2}}{2} \pmod{p}
   \]
when $t_{1} \equiv t_{2} \pmod{p}$, which holds because $p$ is odd, and the relations $\sigma\rme_{i} = \rme_{i+1}$ and $\tau\rme_{i} = \omega^{i}\rme_{i}$). This concludes the proof of the Lemma.
\end{proof}

\begin{corollary}\label{cor:all-conjugate}
Any two elements in $\cycmu$ different from the identity are conjugate in $\pgl$.
\end{corollary}

\begin{remark}
It is interesting to observe that the Proposition, and its Corollary, are false for $p = 2$; then the matrix $\sigma\tau$ has eigenvalues $\pm\sqrt{-1}$, which are not square roots of $1$.
\end{remark}

We will denote by $\xi$ and $\eta$ the first Chern classes in $\ch[1]{\cycmu}$ of the characters $\alpha$ and $\beta$. Then we have
   \[
   \ch{\cycmu} = \ZZ[\xi,\eta]/(p\xi, p\eta).
   \]

We will identify $\cycmu$ with $\FF_{p}\times\FF_{p}$, by sending $\sigma$ to $(1,0)$ and $\tau$ to $(0,1)$; this identifies the automorphism group of $\cycmu$ with $\GL[2](\FF_{p})$.

We are interested in the action of the normalizer $\rN_{\cycmu}\pgl$ of $\cycmu$ in $\pgl$ on $\cycmu$ and on the Chow ring $\ch{\cycmu}$.




\begin{proposition}\label{prop:describe-ch-cycmu}
Consider the homomorphism
   \[
   \rN_{\cycmu}\pgl \arr \GL[2](\FF_{p})
   \]
defined by the action of $\rN_{\cycmu}\pgl$ on $\cycmu$. Its kernel is $\cycmu$, while its image is\/ $\slfinite$.

Furthermore, the ring of invariants
   \[
   \bigl(\ch{\cycmu}\bigr)^{\slfinite}
   \]
is the subring of $\ch{\cycmu}$ generated by the two homogeneous polynomials
   \begin{align*}
      q &\eqdef \eta^{p^{2}-p} + \xi^{p-1}(\xi^{p-1} - \eta^{p-1})^{p-1}\\
      &= \xi^{p^{2}-p} + \eta^{p-1}(\xi^{p-1} - \eta^{p-1})^{p-1}
\intertext{and}
   r &\eqdef \xi\eta(\xi^{p-1} - \eta^{p-1})
   \end{align*}
\end{proposition}

The equality of the two polynomials that appear in the definition of $q$ is not immediately obvious, but is easy to prove, by subtracting them and using the identity
   \[
      (\xi^{p-1} - \eta^{p-1})^{p} =
      \xi^{p^{2}-p} - \eta^{p^{2}-p}.
   \]

\begin{remark}
This is well known to homotopy theorists, as I learnt from N{.} Yagita. More generally, the rings of invariant $\FF_{p}[x_{1}, \dots, x_{n}]^{\GL[n](\FF_{p})}$ and $\FF_{p}[x_{1}, \dots, x_{n}]^{\SL[n](\FF_{p})}$ were computed by L.\,E{.} Dickson, in \cite{dickson-algebra}; the first is known as the \emph{Dickson algebra}. However, I prefer to leave the present self-contained treatment of this case as stands (partly because Lemma~\ref{lem:chern-regular} below will be used several times in what follows).
\end{remark}

\begin{proof} First of all, let us show that the image of the homomorphism above is contained in $\slfinite$. There is canonical symplectic form
   \[
   \bigwedge\nolimits^{2}\bigl(\cycmu\bigr) \arr \mmu
   \]
defined as follows: if $a$ and $b$ are in $\cycmu\subseteq \pgl$, lift them to matrices $\overline{a}$ and $\overline{b}$ in $\GL$. Then the commutator $\overline{a}\overline{b}\overline{a}^{-1}\overline{b}^{-1}$ is a scalar multiple of the identity matrix $\rI_{p}$; it is easy to see that the scalar factor, which we denote by $\generate{a,b}$, is in $\mmu$, and that it only depends on $a$ and $b$, that is, it is independent of the liftings. The resulting function
   \[
   \generate{-,-}\colon (\cycmu)\times(\cycmu) \arr \mmu
   \]
is the desired symplectic form.

Now, $\slfinite$ has $p(p^{2}-1)$ elements. According to Corollary~\ref{cor:all-conjugate}, the action of $\rN_{\cycmu}\pgl$ is transitive on the non-zero vectors in $\FF_{p}^{2}$; so the order of the image of $\rN_{\cycmu}\pgl$ in $\slfinite$ has order divisible by $p^{2}-1$. It is easy to check that the diagonal matrix
   \[
   [1, \omega, \omega^{3}, \dots,
   \underbrace{\omega^{\binom{i}{2}}}_{\text{$i\th$ place}},
   \dots, \omega, 1]
   \]
is also in $\rN_{\cycmu}\pgl$, acts non-trivially on $\cycmu$, and has order $p$. So the order of the image of $\rN_{\cycmu}\pgl$ is divisible by $p$; it follows that it is equal to all of $\slfinite$.

It is not hard to check that the centralizer of $\cycmu$ equals $\cycmu$; and this completes the proof of the first part of the statement.

To study the invariant subring $\bigl(\ch{\cycmu}\bigr)^{\slfinite}$, we use the natural surjective homomorphism
   \[
   \ch{\cycmu} = \ZZ[\xi,\eta]/(p\xi,p\eta) \arr \FF_{p}[\xi,\eta],
   \]
which is an isomorphism in all degrees except $1$; it is enough to show that the  ring of invariants $\FF_{p}[\xi,\eta]^{\slfinite}$ is the polynomial subring $\FF_{p}[q,r]$.

To look for invariants in $\FF_{p}[\xi,\eta]$, we compute the symmetric functions of the vectors in the dual vector space $\bigl(\FF_{p}^{2}\bigr)^{\vee}$; these are the homogeneous components of the polynomial
   \[
   \prod_{i,j \in \FF_{p}}(1 + i\xi + j\eta),
   \]
which are evidently invariant under $\GL[2](\FF_{p})$.

\begin{lemma}\label{lem:chern-regular}\hfil
   \[
      \prod_{0\leq i,j\leq p-1} (1+i\xi+j\eta) = 1 - q + r^{p-1}.
   \]
\end{lemma}

\begin{proof}
Using the formula
   \[
   \prod_{i \in \FF_{p}}(a + ib) = a^{p} - ab^{p-1},
   \]
which holds for any two elements $a$ and $b$ of a commmutative $\FF_{p}$-algebra, we obtain
   \[
      \begin{split}
   \prod_{i,j \in \FF_{p}}(1 + i\xi + j\eta) &=
      \prod_{i \in \FF_{p}}\bigl((1+i\xi)^{p} - (1+i\xi)\eta^{p-1}\bigr)\\
   &= \prod_{i \in \FF_{p}}
      \bigl((1-\eta^{p-1})+i(\xi^{p} -\xi\eta^{p-1})\bigr)\\
   &= 
   (1 - \eta^{p-1})^{p} - (1 - \eta^{p-1})(\xi^{p} -\xi\eta^{p-1})^{p-1}\\
   &= 1 - (\eta^{p^{2}-p} + (\xi^{p} -\xi\eta^{p-1})^{p-1})\\
      &\qquad+ \xi^{p-1}\eta^{p-1}(\xi^{p-1} - \eta^{p-1})^{p-1}\\
   &= 1 - q + r^{p-1}.\qedhere
   \end{split}
   \]
\end{proof}

This shows that $q$ and $r^{p-1}$ are invariant under $\GL[2](\FF_{p})$. The polynomial $r$ is not invariant under $\GL[2](\FF_{p})$, but it is invariant under $\slfinite$. The simplest way to verify this is to observe that $r$ must be a semi-invariant of $\GL[2](\FF_{p})$ (if $g \in \GL[2](\FF_{p})$, then $(gr)^{p-1} = r^{p-1}$, and this means that $gr$ and $r$ differ by a constant in $\FF_{p}^{*}$). But the commutator subgroup of $\GL[2](\FF_{p})$ is well known to be $\slfinite$; so any character $\GL[2](\FF_{p}) \arr \FF_{p}^{*}$ is trivial on $\slfinite$, and $r$ is invariant under $\slfinite$.

We have left to check that $q$ and $r$ generate the ring $\FF_{p}[\xi,\eta]^{\slfinite}$. The equalities
   \[
   \xi^{p^{2}-1} - q\xi^{p-1} + r^{p-1}
   = \eta^{p^{2}-1} - q\eta^{p-1} + r^{p-1}
   = 0,
   \]
which are easily checked by homogenizing the equality of Lemma~\ref{lem:chern-regular}, that is, by adding an indeterminate $t$ and obtaining
   \[
   \prod_{\substack{0\leq i,j\leq p-1\\(i,j)\neq(0,0)}}
   (t+i\xi+j\eta) = t^{p^{2}-1} - qt^{p-1} + r^{p-1},
   \]
ensure that the extension $\FF_{p}[q,r] \subseteq \FF_{p}[\xi,\eta]$ is finite. Hence it is flat, and its degree equals
   \begin{align*}
   \dim_{\FF_{p}}\FF_{p}[\xi,\eta]/(q,r)
   &= \dim_{\FF_{p}}\FF[\xi,\eta]/(q,\xi)\\
      &\qquad + \dim_{\FF_{p}}\FF[\xi,\eta]/(q,\eta)\\
      &\qquad +\dim_{\FF_{p}}\FF[\xi,\eta]
          /\bigl(q,\xi^{p-1} - \eta^{p-1}\bigr)\\
   &= \dim_{\FF_{p}}\FF_{p}[\xi,\eta]/(\eta^{p^{2}-p},\xi)\\
      &\qquad + \dim_{\FF_{p}}\FF_{p}[\xi,\eta]/(\xi^{p^{2}-p},\eta)\\
      &\qquad +\dim_{\FF_{p}}\FF_{p}[\xi,\eta]
          /\bigl(\xi^{p^{2}-p},\xi^{p-1} - \eta^{p-1}\bigr)\\
   &= (p^{2} - p) + (p^{2} - p) + (p^{2} - p)(p - 1)\\
   &= p(p^{2} - 1),
   \end{align*}
which is the order of $\slfinite$. So the degrees of the field extensions
   \[
   \FF_{p}(q,r) \subseteq \FF_{p}(\xi, \eta) 
   \quad\text{and}\quad
   \FF_{p}(\xi,\eta)^{\slfinite} \subseteq \FF_{p}(\xi, \eta)
   \]
both equal $p(p^{2} - 1)$, so $\FF_{p}(q,r) = \FF_{p}(\xi,\eta)^{\slfinite}$; and the result follows, because $\FF_{p}[q,r]$ is integrally closed.
\end{proof}

For later use, let us record the following fact. The image the restriction homomorphism $\ch{\pgl} \arr \ch{\cycmu}$ is contained in $\bigl(\ch{\cycmu}\bigr)^{\slfinite}$. We are going to need formulas for the restrictions of the Chern classes $\chern_{i}(\liesl)$ to $\ch{\cycmu}$.

\begin{lemma}\label{lem:restrict-chern}
Let $i$ be a positive integer. Then the restriction of $\chern_{i}(\liesl)$ to $\ch{\cycmu}$ is $-q$ if $i = p^{2} - p$, is $r^{p-1}$ if $i = p^{2}-1$, and is $0$ in all other cases.
\end{lemma}

\begin{proof}
The total Chern class of $\liegl$ coincides with the total Chern class of $\liesl$, because $\liegl$ is the direct some of $\liesl$ and a trivial representation. From Lemma~\ref{lem:describe-pgl} we see that this total Chern class, when restricted to $\ch{\cycmu}$, equals
   \[
   \sum_{i,j=1}^{p}(1+i\xi+j\eta);
   \]
and then the result follows from Lemma~\ref{lem:chern-regular}.
\end{proof}

We will also need to know about the cohomology ring $\H{\cycmu}$. For any cyclic group $\rC_{n} \simeq \mu_{n}$, the homomorphism $\ch{\rC_{n}} \arr \H{\rC_{n}}$ is an isomorphism. This does not extend to $\cycmu$; however, from the universal coefficients theorem for cohomology, for each index $k$ we have a split exact sequence
   \[
   0\arr \bigoplus_{i+j=k}\H[i]{\cyc}\otimes\H[j]{\mmu}
   \arr \H[k]{\cycmu} \arr
   \bigoplus_{i+j=k+1}\tor_{1}^{\ZZ}\bigl(\H[i]{\cyc}, \H[j]{\mmu}\bigr)
   \arr 0;
   \]
furthermore, since the exterior product homomorphism $\ch{\cyc}\otimes\ch{\mmu} \arr \ch{\cycmu}$ is an isomorphism, the image of the term $\bigoplus_{i+j=k}\H[i]{\cyc}\otimes\H[j]{\mmu}$ into $\H{\cycmu}$ is the image of the cycle homomorphism $\ch{\cycmu} \arr \H{\cycmu}$. From this it is easy to deduce that the cycle homomorphism induces an isomorphism of $\ch{\cycmu}$ with the even dimensional part $\H[\even]{\cycmu}$ of the cohomology.

We have isomorphisms
   \[
   \H[3]{\cycmu} \simeq
   \tor_{1}^{\ZZ}\bigl(\H[2]{\cyc}, \H[2]{\mmu}\bigr)
   \simeq \ZZ/p\ZZ;
   \]
chose a generator $s$ of $\H[3]{\cycmu}$ (later we will make a canonical choice). We have that $s^{2} = 0$, because $p$ is odd, and $s$ has odd degree.

The odd-dimensional part $\H[\odd]{\cycmu}$ of the cohomology is isomorphic to the direct sum $\bigoplus_{i,j}\tor_{1}^{\ZZ}\bigl(\H[i]{\cyc}, \H[j]{\mmu}\bigr)$, with a shift by $1$ in degree. Both $\H[\odd]{\cycmu}$ and $\bigoplus_{i,j} \tor_{1}^{\ZZ}\bigl(\H[i]{\cyc}, \H[j]{\mmu}\bigr)$ have natural structures of modules over $\H{\cyc}\otimes \H{\mmu} = \H[\even]{\cycmu}$, and the isomorphism above is an isomorphism of modules. But $\bigoplus_{i,j}\tor_{1}^{\ZZ}\bigl(\H[i]{\cyc}, \H[j]{\mmu}\bigr)$ is easily seen to be a cyclic $\H[\even]{\cycmu}$-module generated by $s$. From this we obtain the following result.

\begin{proposition}\label{prop:describe-H-cycmu}\hfil
   \[
   \H{\cycmu} = \ZZ[\xi,\eta,s]/(p\xi, p\eta, ps, s^{2}).
   \]
\end{proposition}

We are also interested in the action of $\slfinite$ on $\H{\cycmu}$. I claim that the class $s$ is invariant: this is equivalent to the following.

\begin{lemma}
The action of\/ $\slfinite$ on $\H[3]{\cycmu}$ is trivial.
\end{lemma}

This follows, for example, from the construction of Section~\ref{sec:rho}, where we construct a class $\beta \in \H[3]{\pgl}$ that maps to a non-zero element of $\H[3]{\cycmu}$. It would be logically correct to postpone the proof to Section~\ref{sec:rho}, as this fact is not used before then; but this does not seem very satisfactory, so we prove it now directly.

\begin{proof}
Consider the exact sequence
   \[
   \H[2]{}(\cycmu, \ZZ/p\ZZ) \stackrel{\beta}\longrightarrow
   \H[3]{}(\cycmu, \ZZ) \stackrel{p}\longrightarrow
   \H[3]{}(\cycmu, \ZZ)
   \]
coming from the short exact sequence
   \[
   0 \arr \ZZ \stackrel{p}\longrightarrow \ZZ \arr \ZZ/p\ZZ \arr 0;
   \]
since $\H[3]{}(\cycmu, \ZZ) = \H[3]{\cycmu}$ is $\ZZ/p\ZZ$, we see that the Bockstein homomorphism $\beta\colon \H[2]{}(\cycmu, \ZZ/p\ZZ) \arr \H[3]{}(\cycmu, \ZZ)$ is surjective. It is also $\slfinite$-equivariant. By K\"unneth's formula, the exterior product induces an isomorphism of the direct sum
   \[
   \H[2]{}(\cyc, \ZZ/p\ZZ)\\ \oplus
   \bigl(\H[1]{}(\cyc, \ZZ/p\ZZ) \otimes \H[1]{}(\mmu, \ZZ/p\ZZ)\bigr)
   \oplus 
   \H[2]{}(\mmu, \ZZ/p\ZZ)
   \]
with $\H[2]{}(\cycmu, \ZZ/p\ZZ)$. Now, from the commutativity of the diagram
   \[
   \xymatrix{
      \H[2]{}(\cyc, \ZZ/p\ZZ)\ar[r]\ar[d] & \H[3]{}(\cyc, \ZZ)\ar[d]
      \hsmash{{}= 0}\\
   \H[2]{}(\cycmu, \ZZ/p\ZZ)\ar[r]^-{\beta} & \H[3]{}(\cycmu, \ZZ)
   }
   \]
where the rows are Bockstein homomorphisms and the columns are induced by projection $\cycmu \arr \cyc$, we see that the Bockstein homomorphism $\beta$ sends $\H[2]{}(\cyc, \ZZ/p\ZZ)$, and $\H[2]{}(\mmu, \ZZ/p\ZZ)$ for analogous reasons, to $0$. Hence the composite of the exterior product map
   \[
   \H[1]{}(\cyc, \ZZ/p\ZZ) \otimes \H[1]{}(\mmu, \ZZ/p\ZZ) \arr
      \H[2]{}(\cycmu, \ZZ/p\ZZ)
   \]
with $\beta$ is surjective. But we have an isomorphism
   \[
   \H[1]{}(\cycmu, \ZZ/p\ZZ) \simeq 
   \H[1]{}(\cyc, \ZZ/p\ZZ) \oplus \H[1]{}(\mmu, \ZZ/p\ZZ)
   \]
which induces an isomorphism
   \[
   \bigwedge\nolimits^{2}\H[1]{}(\cycmu, \ZZ/p\ZZ) \simeq 
   \H[1]{}(\cyc, \ZZ/p\ZZ) \otimes \H[1]{}(\mmu, \ZZ/p\ZZ).
   \]
This shows that the composite of the map
   \[
   \bigwedge\nolimits^{2}\H[1]{}(\cycmu, \ZZ/p\ZZ) \arr \H[2]{}(\cycmu,\ZZ/p\ZZ)
   \]
with the Bockstein homomorphism $\beta$ is surjective, hence an isomorphism, because both groups are isomorphic to $\ZZ/p\ZZ$. It is also evidently $\GL[2](\FF_{2})$-equivariant. The action of $\GL[2](\FF_{p})$ on $\bigwedge\nolimits^{2}\H[1]{}(\cycmu, \ZZ/p\ZZ)$ is by multiplication by the inverse of the determinant; hence $\slfinite$ acts trivially, and this completes the proof.
\end{proof}

From this we deduce the following fact.

\begin{proposition}\label{prop:H-invariants-cycmu}
The ring of invariants $\bigl(\H{\cycmu}\bigr)^{\slfinite}$ is generated by $q$, $r$ and $s$.
\end{proposition}

\begin{remark}
The group $\cycmu$ is important in the theory of division algebras. Suppose that $K$ is a field containing $k$, and $E \arr \spec K$ is a non-trivial $\pgl$ principal bundle. This corresponds to a central division algebra $D$ over $K$ of degree $p$. Recall that $D$ is \emph{cyclic} when there are elements $a$ and $b$ of $K^{*}$, such that $D$ is generated by two elements $x$ and $y$, satisfying the relations $x^{p} = a$, $y^{p} = b$, $yx = \omega xy$. It is not hard to show that $D$ is cyclic if and only if $E$ has a reduction of structure group to $\cycmu$.

One of the main open problems in the theory of division algebra is whether all division algebras of prime degree are cyclic. Let $V$ be a representation of $\pgl$ over $k$ with a non-empty open invariant subset $U$ on which $\pgl$ acts freely. Let $K$ be the fraction field of $U/\pgl$, $E$ the pullback to $\spec K$ of the $\pgl$-torsor $U \arr U/G$ and $D$ the corresponding division algebra; it is well known that $D$ cyclic if and only if every division algebra of degree $p$ over a field containing $k$ is cyclic. 

The obvious way to show that $D$ is \emph{not} cyclic is to show that there is an invariant for division algebras that is $0$ for cyclic algebras, but not $0$ for $D$. However, the result proved here implies that there is no such invariant in the cohomology ring $\H{\pgl}$. In fact, consider a non-zero invariant $\xi \in\H{\pgl}$. Then either $\xi$ has even degree, so it comes from $\ch{\pgl}$, hence it restricts to $0$ in $V/\pgl$ for some open invariant subset $V \subseteq U$, or it has odd degree, and then it maps to $0$ in $\ch{\toruspgl}$, and it does not map to $0$ in $\ch{\cycmu}$.

This is related with the fact that one can not find such an invariant in \'etale cohomology with $\ZZ/p\ZZ$ coefficients (see \cite[\S 22.10]{garibaldi-merkurjev-serre}).
\end{remark}

\section{On $\cycgl$}\label{sec:cycgl}

\begin{proposition}\label{prop:chcycgl-hcycgl}
Assume that $k = \CC$. Then the cycle homomorphism $\ch{\cycgl} \to \H{\cycgl}$ is an isomorphism.
\end{proposition}

\begin{proof} This the first illustration of the stratification method: we take a geometrically meaningful representation of $\cycgl$ and we stratify it.

Denote by $V \eqdef \AA^p$ the standard representation of $\GL$, restricted to $\cycgl$. We denote by $V_{\le i}$ the Zariski open  $\cycgl$-invariant subset consisting of $p$-uples of complex numbers such that at most $i$ of them are $0$, and by $V_i \eqdef V_{\le i} \setminus V_{\le i-1}$ the smooth locally closed subvariety of $p$-uples consisting of vectors with exactly $i$ coordinates that are $0$. Obviously $V_{\le p-1} = V \setminus \{0\}$ and $V_p = 0$.

\begin{lemma}\label{lem:aux1-isom}
For each $0 \le i \le p-1$, the cycle homomorphism $\ch{\cycgl}{V_i} \to \H{\cycgl}{V_i}$ is an isomorphism.
\end{lemma}

\begin{proof}
First of all, assume that $i = 0$. Then the action of $\cycgl$ on $V_0$ is transitive, and the stabilizer of $(1, \dots, 1) \in V_{0}(k)$ is $\cyc$; hence  we have a commutative diagram
   \[
   \xymatrix{
   \ch{\cycgl}(V_0)\ar[r]\ar[d] & \H{\cycgl}(V_0)\ar[d]\\
   \ch{\cyc}         \ar[r]       &  \H{\cyc}
   }
   \]
where the rows are cycle homomorphisms and the columns are isomorphisms. Since the bottom row is also an isomorphism, the thesis follows.

When $i>0$ the argument is similar. The action of $\cycgl$ on $V_i$ expresses $V_i$ as a disjoint union of open orbits $\Omega_1$, \dots,~$\Omega_r$, where $r \eqdef \frac{1}{p} \binom{p}{i}$, and the stabilizer of a point of each $\Omega_j$ is an $i$-dimensional torus $T_j$; hence we get a commutative diagram
   \[
   \xymatrix{
   \ch{\cycgl}(V_i)\ar[r]    \ar[d] & \H{\cycgl}(V_i)\ar[d]\\
   \bigoplus_{h=1}^r\ch{T_j} \ar[r] & \bigoplus_{h=1}^r\H{T_j}
   }
   \]
where the columns and the bottom row are isomorphisms.
\end{proof}

\begin{lemma}\label{lem:aux2-isom}
For each $0 \le i \le p-1$, the cycle homomorphism $\ch{\cycgl}{V_{\le i}} \to \H{\cycgl}{V_{\le i}}$ is an isomorphism.
\end{lemma}

\begin{proof}

We proceed by induction on $i$. When $i = 0$ we have $V_{\le 0} = V_0$, and the thesis follows from the previous lemma. For the inductive step, we have a commutative diagram with exact rows
   \[
   \xymatrix{
   \ch{\cycgl}(V_i)         \ar[r]\ar[d]             &
   \ch{\cycgl}(V_{\le i})   \ar[r]\ar[d]             &
   \ch{\cycgl}(V_{\le i-1}) \ar[r]\ar[d]^*+[o][F]{1} & 0 \\
   \H{\cycgl}(V_i)         \ar[r]  ^*+[o][F]{3}      &
   \H{\cycgl}(V_{\le i})   \ar[r]^*+[o][F]{2}        &
   \H{\cycgl}(V_{\le i-1}) \hsmash{;}
   }
   \]
by inductive hypothesis, the arrow marked with $\xymatrix{*+[o][F]{1}}$ is an isomorphism, hence the arrow marked with $\xymatrix{*+[o][F]{2}}$ is surjective. However, the bottom row of the diagram extends to a Gysin exact sequence
   \[
   \H{\cycgl}(V_i) \to
   \H{\cycgl}(V_{\le i}) \to
   \H{\cycgl}(V_{\le i-1}) \to
   \H{\cycgl}(V_i) \to \cdots
   \]
showing that the arrow marked with $\xymatrix{*+[o][F]{3}}$ is injective. From this, and the fact that the left hand column is an isomorphism, it follows that the middle column is also an isomorphism, as desidered.
\end{proof}

Let us proceed with the proof of the Theorem. For each $i$ we have an commutative diagram with exact rows
   \[
   \xymatrix{
   \ch{\cycgl} \ar[r]^{\c_p (V)}\ar[r]\ar[d] &
   \ch{\cycgl} \ar[r]\ar[r]\ar[d]&
   \ch{\cycgl} (V\setminus \{0\}) \ar[r]\ar[r]\ar[d]&
   0\\
   \H{\cycgl} \ar[r]^{\c_p (V)}\ar[r] &
   \H{\cycgl} \ar[r]\ar[r]&
   \H{\cycgl} (V\setminus \{0\})\\
   }
   \]
Now, by Lemma~\ref{lem:aux2-isom} the right hand column is an isomorphism, hence, arguing as in the proof of Lemma~\ref{lem:aux2-isom}, we conclude that the bottom row of the diagram is a short exact sequence.

If $i$ is odd, we have $\H[i]{\cycgl} (V\setminus \{0\}) = 0$, hence the multiplication homomorphism 
   \[
   \H[i-2p]{\cycgl} \xarr{\c_p (V)} \H[i]{\cycgl}
   \]
is an isomorphism. From this we deduce, by induction on $i$, that $\H[i]{\cycgl} = 0$ for all odd $i$.

When $i$ is even, one proceeds similarly by induction on $i$, with a straightforward diagram chasing in the  diagram above.
\end{proof}

Let us compute the Chow ring of the classifying space of $\cycgl$. The Weyl group $\s$ acts on $\ch{\torusgl} = \ZZ[x_{1}, \dots, x_{p}]$ by permuting the $x_i$'s. Consider the action of $\cyc$ on $\ch{\torusgl}$: the group permutes the monomials, and the only monomials that are left fixed are the ones of the form $\sigma_{p}^{r} = x_1^r \dots x_p^r$, while on the others the action of $\cyc$ is free. We will call the monomials that are not powers of $\sigma_{p}$ \emph{free monomials}. Then $\ch{\torusgl}$ splits as a direct sum $\ZZ[\sigma_{p}]\oplus M$, where $M$ is the free $\ZZ\cyc$-module generated by the free monomials. Hence the ring of invariants $(\ch{\torusgl})^{\cyc}$ is a direct sum $\ZZ[\sigma_{p}] \oplus M^{\cyc}$, and $M^{\cyc}$ is a free abelian group on the generators $\sum_{s \in \cyc} sm$, where $m$ is a free monomial.

We will denote by $\xi \in \ch[1]{\cyc}$ the first Chern class of the character $\cyc \to \gm$ obtained by sending the generator $(1, \dots, p)$ of $\cyc$ to the fixed generator $\omega$ of $\mmu$, and also its pullback to $\cycgl$ through the projection $\cycgl \to \cyc$.

We will also use the subgroup $\mmu \subseteq \torusgl$ of matrices of the form $\zeta\rI_{p}$, where $\zeta\in \mmu$. The Chow ring $\ch{\mmu}$ is of the form $\ZZ[\eta]/(p\eta)$, where $\eta$ is the first Chern class of the 1-dimensional representation given by the embedding $\mmu \into \gm$. The action of $\cyc$ on $\mmu$ is trivial, so there a copy of $\cycmu$ in $\cycgl$; the Chow ring $\ch{\cycmu}$ is $\ZZ[\xi, \eta]/(p\xi, p\eta)$.

Here are the facts about $\ch{\cycgl}$ that we are going to need.

\begin{proposition}\call{prop:describe-cycgl}\hfil
\begin{enumeratea}

\itemref{4} The image of the restriction homomorphism
   \[
   (\ch{\torusgl})^{\cyc} \arr \ch{\mmu}  = \ZZ[\eta]/(p\eta)
   \]
is the subring generated by $\eta^{p}$, which is the image of $\sigma_{p}$. The kernel is the subgroup of $(\ch{\torusgl})^{\cyc}$ generated by the $\sum_{s \in \cyc} sm$, where $m$ is a free monomial, and by $p\sigma_{p}$.

\itemref{1} The ring homomorphism $\ch{\cycgl} \arr \bigl(\ch{\torusgl}\bigr)^{\cyc}$ induced by the embedding $\torusgl \into \cycgl$ is surjective, and admits a canonical splitting $\phi\colon \bigl(\ch{\torusgl}\bigr)^{\cyc} \arr \ch{\cycgl}$, which is a ring homomorphism.

\itemref{2} As an algebra over $\ch{\cycgl}$, the ring $\ch{\cycgl}$ is generated by the element $\xi$, while the ideal of relations is generated by the following: $p\xi = 0$, and $\phi(u) \xi = 0$ for all $u$ in the kernel of the ring homomorphism $\ch{\torusgl} \arr \ch{\mmu}$ induced by the embedding $\mmu \into \torusgl$.

\itemref{3} The ring homomorphism $\ch{\cycgl} \arr \ch{\torusgl} \times \ch{\cycmu}$ induced by the embeddings $\torusgl \into \cycgl$ and $\cycmu \into \cycgl$ is injective.

\itemref{5} The restriction homomorphism $\ch{\cycgl} \arr (\ch{\torusgl})^{\cyc}$ sends the kernel of $\ch{\cycgl} \arr \ch{\cycmu}$ bijectively onto the kernel of $\ch{\torusgl} \arr \ch{\mmu}$.

\end{enumeratea}
\end{proposition}

\begin{proof}
Let us prove part~\refpart{prop:describe-cycgl}{4}.
All the $x_{i}$ in $\ch{\torusgl}$ map to $\eta$ in $\ch{\mmu}$, so $\sigma_{p}$ maps to $\eta^{p}$, and all the $\sum_{s \in \cyc} sm$ map to $p\eta^{\deg m} = 0$.

Let us prove \refpart{prop:describe-cycgl}{1}. First of all let us construct the splitting $\phi\colon \bigl(\ch{\torusgl}\bigr)^{\cyc} \arr \ch{\cycgl}$ as a homomorphism of abelian groups. The group $\bigl(\ch{\torusgl}\bigr)^{\cyc}$ is free over the powers of $\sigma_{p}$ and the $\sum_{s \in \cyc}sm$.

The restriction of the canonical representation $V$ of $\cycgl$ to the maximal torus $\torusgl$ splits are a direct sum of 1-dimensional representations with first Chern characters $x_{1}$, \dots,~$x_{p}$; hence the $i\th$ Chern class $\c_{i}(V) \in \ch[i]{\cycgl}$ restricts to $\sigma_{i} \in (\ch[i]{\torusgl})$. We define the splitting by the rules
\begin{enumeratea}

\item $\phi(\sigma_{p}^{r}) = \c_{p}(V)^{r} \in \ch{\cycgl}$ for each $r > 0$, and

\item $\phi(\sum_{s \in \cyc}sm) = \tsf{\torusgl}{\cycgl} m \in \ch{\cycgl}$
for each free monomial $m$.

\end{enumeratea}

Notice that the transfer in the second part of the definition only depends on the orbit of $m$; hence $\phi$ is well defined.

We need to check that $\phi$ is a ring homomorphism, by taking two basis element $u$ and $v$ and showing that $\phi(uv)$ equals $\phi(u)\phi(v)$. This is clear when both $u$ and $v$ are powers of $\sigma_{p}$.

Consider the product $\sigma_{p}^{r}\sum_{s\in \cyc}sm = \sum_{s\in \cyc}s(\sigma_{p}^{r}m)$; we have
   \begin{align*}
   \phi\Bigl(\sigma_{p}^{r}\sum_{s\in \cyc}sm\Bigr) &= 
      \phi\Bigl(\sum_{s\in \cyc}s(\sigma_{p}^{r}m)\Bigr)\\
   & = \tsf{\torusgl}{\cycgl} (\sigma_{p}^{r}m)
      \quad\text{(because $\sigma_{p}^{r}m$ is still a free monomial)}\\
   & = \c_{p}(V)^{r}\tsf{\torusgl}{\cycgl} (m)
      \quad\text{(by the projection formula)}\\
   & = \phi(\sigma_{p}^{r})\phi\Bigl(\sum_{s\in \cyc}sm\Bigr).
   \end{align*}

Now the hardest case. Notice that if $m$ is any monomial, not necessarily free, we have the equality
   \[
   \phi\Bigl(\sum_{s\in \cyc}sm\Bigr) = \tsf{\torusgl}{\cycgl} m.
   \]
When $m$ is free this holds by definition, whereas when $m = \sigma_{p}^{r}$ we have
   \begin{align*}
   \smash{\phi\Bigl(\sum_{s\in \cyc}\sigma_{p}^{r}\Bigr)}
   &= p \phi(\sigma_{p}^{r})\\
   &= p\c_{p}(V)^{r} \\
   &= \tsf{\torusgl}{\cycgl} \res{\cycgl}{\torusgl}\c_{p}(V)^{r}\\
   &= \tsf{\torusgl}{\cycgl} \sigma_{p}^{r}.
   \end{align*}
Take two free monomials $m$ and $n$. We have
   \begin{align*}
   \phi\Bigl(\sum_{s\in \cyc}sm \cdot\sum_{s\in \cyc}sn\Bigr)
      &= \phi\Bigl(\sum_{s,t \in \cyc}sm \cdot tn\Bigr)\\
   &= \phi\Bigl(\sum_{s,t \in \cyc}sm \cdot stn\Bigr)\\
   &= \sum_{t\in \cyc}\phi\Bigl(\sum_{s \in \cyc}s(m \cdot tn)\Bigr)\\
   &= \sum_{t \in \cyc}\tsf{\torusgl}{\cycgl}(m\cdot tn)\\
   &= \tsf{\torusgl}{\cycgl}\Bigl(m\cdot \sum_{t \in \cyc}tn\Bigr)\\
   &= \tsf{\torusgl}{\cycgl}
      \Bigl(m\cdot \res{\cycgl}{\torusgl}\tsf{\torusgl}{\cycgl}n\Bigr)\\
   &= \tsf{\torusgl}{\cycgl}(m) \tsf{\torusgl}{\cycgl}(n)\\
   &= \phi\Bigl(\sum_{s\in \cyc}sm\Bigr) \phi\Bigl(\sum_{s\in \cyc}sn\Bigr)
   \end{align*}
as claimed. This ends the proof of part~\refpart{prop:describe-cycgl}{1}.

For parts~\refpart{prop:describe-cycgl}{2} and \refpart{prop:describe-cycgl}{3}, notice the following fact: since the restriction of $\xi$ to $\ch{\torusgl}$ is $0$, from the projection formula it follows that $\xi \tsf{\torusgl}{\cycgl}(m) = 0 \in \ch{\cycgl}$ for any $m \in \ch{\torusgl}$; hence we get that $\phi(\sum_{s \in \cyc}m)\xi = 0 \in \ch{\cycgl}$, as claimed. Thus, the relations of the statement of the Proposition hold true.

Denote by $\ch[+]{\torusgl}$ the ideal of $\ch{\torusgl}$ generated by homogeneous elements of positive degree. Then the image of $\bigl(\ch[+]{\torusgl}\bigr)^{\cyc}$ in $\ch{\cycgl}$ via $\phi$ maps to $0$ under the restriction homomorphism $\res{\cycgl}{\cyc} \colon \ch{\cycgl} \arr \ch{\cyc}$. In fact, the image of $\ch[+]{\torusgl}$ is generated by elements of the form $\tsf{\torusgl}{\cycgl}m$, where $m \in \ch{\torusgl}$ is a monomial of positive degree, and by positive powers $\c_{p}(V)^{r}$ of the top Chern class of $V$. The fact that the restriction of $\tsf{\torusgl}{\cycgl}m$ is $0$ follows from Mackey's formula. On the other hand, the restriction of $V$ to $\cyc$ is a direct sum of 1-dimensional representations with first Chern classes $0$, $\xi$, $2\xi$, \dots,~$(p-1)\xi$, so the restriction of $V$ has trivial top Chern class.

\begin{lemma}\label{lem:keyfact-cycgl}
The kernel of the restriction homomorphism
   \[
   \res{\cycgl}{\cyc} \colon \ch{\cycgl} \arr \ch{\cyc}
   \]
consists of the sum the image of $(\ch[+]{\torusgl})^{\cyc}$ in $\ch{\cycgl}$ via $\phi$, and of the ideal $\bigl(\c_{p}(V)\bigr) \subseteq \ch{\cycgl}$.
\end{lemma}

\begin{proof}
Consider the hyperplane $H_{i}$ in the canonical representation $V = \AA^{p}$ defined by the vanishing of the $i\th$ coordinate. Denote by $H = \cup_{i=1}^{p} H_{i} \subseteq V$ the union. If $V_{0} = V \setminus H$ we have an exact sequence
   \[
   \ch{\cycgl}(H)  \arr \ch{\cycgl}(V) \arr \ch{\cycgl}(V_{0}) \arr 0.
   \]
We identify $\ch{\cycgl}(V)$ with $\ch{\cycgl}$ via the pullback $\ch{\cycgl} \arr \ch{\cycgl}(V)$, which is an isomorphism. The action of $\torusgl$ on $V_{0}$ is free and transitive, and the stabilizer of the point $(1, 1, \dots, 1)$ is $\cyc \subseteq \cycgl$. Hence we have an isomorphism of $\ch{\cycgl}(V_{0})$ with $\ch{\cyc}$, and the pullback $\ch{\cycgl}(V) \arr \ch{\cycgl}(V_{0})$ is identified with the restriction homomorphism $ \ch{\cycgl} \arr \ch{\cyc}$. So the kernel of this restriction is the image of $\ch{\cycgl}(H)$.

Denote by $\widetilde{H}$ the disjoint union $\coprod_{i=1}^{p}H_{i} \sqcup \{0\}$ of the $H_{i}$ with the origin $\{0\}\subseteq V$. I claim that the proper pushforward $\ch{\cycgl}(\widetilde{H}) \arr \ch{\cycgl}(H)$ is surjective. This follows from Proposition~\ref{prop:envelope}: we need to check that every $\cycgl$-invariant closed subvariety of $H$ is the birational image of a $\cycgl$-invariant subvariety of $\widetilde{H}$. Denote by $W$ a  $\cycgl$-invariant closed subvariety of $H$. If $W = \{0\}$ we are done. Otherwise it is easy to see that $W$ will be the union of $p$ $\torusgl$-invariant irreducible components $W_{1}$, \dots,~$W_{p}$, such that each $W_{i}$ is contained in $H_{i}$. Then the disjoint union $\coprod_{i=1}^{p}W_{i} \subseteq \coprod_{i=1}^{p}H_{i} \subseteq \widetilde{H}$ is $\cycgl$-invariant and maps birationally onto $W$. Hence we conclude that the kernel of the restriction homomorphism is the sum of the images of the proper pushforwards    \[
   \ch{\cycgl}(\{0\}) \arr \ch{\cycgl}(V)
   \]
and
   \[
   \ch{\cycgl}\bigl(\coprod_{i=1}^{p}H_{i}\bigr) \arr \ch{\cycgl}(V).
   \]

After identifying $\ch{\cycgl}(V)$ with $\ch{\cycgl}$, the first pushforward is just multiplication by $\c_{p}(V)$, so its image is the ideal $\bigl(\c_{p}(V)\bigr) \subseteq \ch{\cycgl}$.

Notice that the disjoint union $\coprod_{i=1}^{p} H_{i}$ is canonically isomorphic, as a $\cycgl$-scheme, to $(\cycgl)\times^{\torusgl}H_{1}$; hence there is a canonical isomorphism
    \[
   \ch{\cycgl}\biggl(\,\coprod_{i=1}^{p}H_{i}\biggr) \simeq 
   \ch{\torusgl}(H_{1}).
   \]
The pushforward $\ch{\torusgl}(H_{1}) \arr \ch{\cycgl}(V)$ is the composite of the proper pushforward $\ch{\torusgl}(H_{1}) \arr \ch{\torusgl}(V)$, followed by the transfer homomorphism $\ch{\torusgl}(V) \arr \ch{\cycgl}(V)$. After identifying $\ch{\torusgl}(H_{1})$ and $\ch{\torusgl}(V)$ with $\ch{\torusgl}$, $\ch{\cycgl}(V)$ with $\ch{\cycgl}$, we see that this implies that the image of $\ch{\cycgl}\bigl(\coprod_{i=1}^{p}H_{i}\bigr)$ in $\ch{\cycgl}(V) = \ch{\cycgl}$ is the image of the ideal $(x_{1}) \subseteq \ch{\torusgl}$ under the transfer map $\ch{\torusgl} \arr \ch{\cycgl}$. So each element of the image of $\ch{\cycgl}\bigl(\coprod_{i=1}^{p}H_{i}\bigr)$ can be written as a linear combination with integer coefficients of transfers of monomials of positive degree: and this completes the proof of the Lemma.
\end{proof}

Now we show that $\ch{\cycgl}$ is generated, as an algebra over $\bigl(\ch{\torusgl}\bigr)^{\cyc}$, by the single element $\xi$. Take an element $\alpha$ of $\ch{\cycgl}$ of degree $d$, and write its image in $\ch{\cyc} = \ZZ[\xi]/(p\xi)$ in the form $m\xi^{d}$, where $m$ is an integer. Then $\alpha - m\xi^{d} \in \ch{\cycgl}$ maps to $0$ in $\ch{\cyc}$, so according to Lemma~\ref{lem:keyfact-cycgl} it is of the form $\beta + \sigma_{p}\gamma$, where $\beta$ is in $\bigl(\ch{\torusgl}\bigr)^{\cyc}$ and $\gamma \in \ch[d-p]{\cyc}$. The proof is concluded by induction on $d$.

Now we prove that the relations indicated generate the ideal of relations, and, simultaneously, part~\refpart{prop:describe-cycgl}{3}.

Take an element $\alpha \in \ch[d]{\cycgl}$; using the given relations, we can write $\alpha$ in the form $\alpha_{0} + \alpha_{1}\xi + \alpha_{2}\xi^{2} + \cdots$, where $\alpha_{0} \in \bigl(\ch[d]{\torusgl}\bigr)^{\cyc}$, while for each $i > 0$ the element $\alpha_{i}$ is of the form $d_{i}\sigma_{p}^{r}$, where $0 \leq d_{i} \leq p-1$, and $rp = d-i$, when $p$ divides $d-i$, and $0$ otherwise.

Assume that the image of $\alpha$ in $\ch{\torusgl} \times \ch{\cycmu}$ is $0$. The image of $\alpha$ in $\ch{\torusgl}$ is $\alpha_{0}$, hence $\alpha_{0} = 0$. 

\begin{lemma}\label{lem:restrictcanonical-mmu}
The restriction of $\phi(\sigma_{p}) = \c_{p}(V)$ to $\ch{\cycmu} = \ZZ[\xi, \eta](p\xi, p\eta)$ equals $\eta^{p} - \eta\xi^{p-1}$.
\end{lemma}

\begin{proof}
The restriction of $V$ to $\cycmu$ decomposes as a direct sum of 1-dimensional representations with first Chern classes $\eta$, $\eta - \xi$, $\eta - 2\xi$, \dots,~$\eta - (p-1)\xi$, and
   \[
   \eta (\eta - \xi)(\eta - 2\xi) \dots \bigl(\eta - (p-1)\xi\bigr) = 
   \eta^{p} - \eta\xi^{p-1}. \qedhere
   \]
\end{proof}

Since $\xi$ and $\eta^{p} - \eta\xi^{p-1}$ are algebraically independent in the polynomial ring $\FF_{p}[\xi, \eta]$, it follows that all the $\alpha_{i}$ are all $0$. This finishes the proof of \refpart{prop:describe-cycgl}{2} and \refpart{prop:describe-cycgl}{3}.

Finally, let us prove part~\refpart{prop:describe-cycgl}{5}.

Injectivity follows immediately from part~\refpart{prop:describe-cycgl}{3}. To show that the restriction homomorphism is surjective, it is sufficient to show that if $u$ is in the kernel of the homomorphism $\bigl(\ch{\torusgl}\bigr)^{\cyc} \arr \ch{\mmu}$, then $\phi(u)$ is in the kernel of $\ch{\cycgl} \arr \ch{\cycmu}$. Each element of $\bigl(\ch{\torusgl}\bigr)^{\cyc}$ of the form $\sum_{s \in \cyc}sm$ goes to $0$ in $\ch{\cyc}$, while $\sigma_{p}$ goes to $\eta^{p}$; hence $u$ is a linear combination of elements of the form $\sum_{s \in \cyc}m$ and $p\sigma_{p}^{r}$. So $\phi(u)$ is a linear combination of element of $\ch{\cycgl}$ of the form $p\chern_{p}(V)^{r}$ and $\tsf{\torusgl}{\cycgl}m$; from the following Lemma we see that all these elements to $\ch{\cycmu}$ is $0$.

\begin{lemma}\label{lem:restriction-zero}
If $u$ is an element of positive degree in $\ch{\torusgl}$, the restriction of $\tsf{\torusgl}{\cycgl}u$ to $\ch{\cycmu}$ is $0$.
\end{lemma}
\begin{proof}
The double coset space $(\cycmu)\backslash (\cycgl) /\torusgl$ consists of a single point and $(\cycmu) \cap \torusgl = \mmu$, so we have
   \[
   \res{\cycgl}{\cycmu} \tsf{\torusgl}{\cycgl}u 
   = \tsf{\mmu}{\cycmu}\res{\torusgl}{\mmu} u.
   \]
However, I claim that the transfer homomorphism
   \[
   \tsf{\mmu}{\cycmu}\colon \ch{\mmu} \arr \ch{\cycmu}
   \]
is $0$ in positive degree. In fact, the restriction homomorphism
   \[
   \res{\cycmu}{\mmu}\colon \ch{\cycmu} \arr \ch{\mmu}
   \]
is surjective, because the embedding $\mmu \into \cycmu$ is split by the projection $\cycmu \arr \mmu$. It follows immediately, again from Mackey's formula, that the composition $\tsf{\mmu}{\cycmu}\res{\cycmu}{\mmu}$ is multiplication by $p$; and all classes in $\ch{\cycmu}$ in positive degree are $p$-torsion.
\end{proof}

This concludes the proof of Proposition~\ref{prop:describe-cycgl}.
\end{proof}

\begin{remark}
When $k = \CC$, Propositions \ref{prop:chcycgl-hcycgl} and \ref{prop:describe-cycgl} give a description of the cohomology $\H{\cycgl}$. This can be proved directly, by studying the Hochschild--Serre spectral sequence
   \[
   E_{2}^{ij} = \rH^{i}\bigl(\cyc, \rH^{j}_{\torusgl}\bigr) \Longrightarrow
      \H[i+j]{\cycgl}.
   \]
\end{remark}

\section{On $\cycpgl$}\label{sec:cycpgl}

In this section we study the Chow ring of the classifying space of the group $\cycpgl$. Here is our main result. Consider the subgroup $\mmu \subseteq \toruspgl$ defined, as in the Introduction, by the formula $\zeta \mapsto [\zeta, \zeta^{2}, \dots, \zeta^{p-1}, 1]$. This defines a homomorphism of rings $\ch{\toruspgl} \arr \ch{\mmu}$. 

\begin{proposition}\call{prop:describe-cycpgl}\hfil
\begin{enumeratea}

\itemref{1} The image of the restriction homomorphism
   \[
   (\ch{\toruspgl})^{\cyc} \arr \ch{\mmu}  = \ZZ[\eta]/(p\eta)
   \]
is the subring generated by $\eta^{p}$.

\itemref{2} The ring homomorphism $\ch{\cycpgl} \arr \bigl(\ch{\toruspgl}\bigr)^{\cyc}$ induced by the embedding $\toruspgl \into \cycpgl$ is surjective, and admits a canonical splitting $\phi\colon \bigl(\ch{\toruspgl}\bigr)^{\cyc} \arr \ch{\cycpgl}$, which is a ring homomorphism.

\itemref{3} As an algebra over $\bigl(\ch{\toruspgl}\bigr)^{\cyc}$, the ring $\ch{\cycpgl}$ is generated by the element $\xi$, while the ideal of relations is generated by the following: $p\xi = 0$, and $\phi(u) \xi = 0$ for all $u$ in the kernel of the ring homomorphism $\ch{\toruspgl} \arr \ch{\mmu}$ induced by the embedding $\mmu \into \toruspgl$.

\itemref{4} The ring homomorphisms
   \[
   \ch{\cycpgl} \arr \ch{\toruspgl} \times \ch{\cycmu}
   \]
and
   \[
   \H{\cycpgl} \arr \H{\toruspgl} \times \H{\cycmu}
   \]
induced by the embeddings $\toruspgl \into \cycpgl$ and $\cycmu \into \cycpgl$ is injective.

\itemref{5} The restriction homomorphism $\ch{\cycpgl} \arr (\ch{\toruspgl})^{\cyc}$ sends the kernel of $\ch{\cycpgl} \arr \ch{\cycmu}$ bijectively onto the kernel of $\ch{\toruspgl} \arr \ch{\mmu}$.
\end{enumeratea}

\end{proposition}

\begin{proof}
One of the main ideas in the paper is to exploit the fact, already used in \cite{vezzosi-pgl3} and rediscovered in \cite{vavpetic-viruel}, that there is an isomorphism of tori
   \[
   \Phi \colon \toruspgl \simeq \torussl
   \]
defined by
   \[
   \Phi(t_1, \dots , t_p) = [t_1/t_p, t_2/t_1, t_2/t_2, \dots,
   t_{p-1}/t_{p-2}, t_p/t_{p-1}].
   \]
This isomorphism is not $\s$-equivariant, but it is $\cyc$-equivariant; therefore it induces an isomorphism
   \[
   \Phi \colon \cycpgl \simeq \cycsl.
   \]

The composite of the embedding $\mmu \into \toruspgl$ with the isomorphism $\Phi$ is the embedding $\mmu \into \torussl$ defined by $\zeta \mapsto [\zeta, \zeta, \dots, \zeta]$.

Now, take an open subset $U$ of a representation of $\cycgl$ on which $\cycgl$ acts freely. The projection $U/\cycsl \arr U/\cycgl$ is a $\gm$-torsor, coming from the determinant $\det \colon \cycgl \arr \gm$ of the canonical representation $V$ of $\torusgl$. Lemma~\ref{lem:exact-top} implies that there is an exact sequence
   \[
   \ch{\cycgl} \xarr{\c_{1}(V)} \ch{\cycgl} \arr \ch{\cycsl} \arr 0
   \]
and a ring isomorphism $\ch{\cycsl} \simeq \ch{\cycgl}/\bigl(\c_{1}(V)\bigr)$.

Consider the splitting $\phi \colon \bigl(\ch{\torusgl}\bigr)^{\cyc} \arr \ch{\cycgl}$ constructed in the previous section. I claim that $\c_{1}(V)$ coincides with $\phi(\sigma_{1}) = \tsf{\torusgl}{\cycgl}x_{1}$. To prove this it is enough, according to Proposition \refall{prop:describe-cycgl}{3}, to show that these two classes coincide after restriction to $\ch{\torusgl}$ and to $\ch{\cycmu}$. The restrictions of both classes to $\ch{\torusgl}$ coincide with $x_{1} + \dots + x_{p}$.

The action of $\cycmu$ on $V$ splits as a direct sum of 1-dimensional representations with characters $\eta + \xi$, $\eta + 2\xi$, \dots,~$\eta + (p-1)\xi$, $\eta$, so the restriction of $\c_{1}(V)$ to $\ch{\cycmu}$ is
   \[
   \eta + \xi + \eta + 2\xi + \dots + \eta + (p-1)\xi + \eta = 
       p\eta + \frac{p(p-1)}{2}\xi = 0.
   \]
So we need to show that the restriction of $\tsf{\torusgl}{\cycgl}x_{1}$ to $\ch{\cycmu}$ is also $0$. This is a particular case of Lemma~\ref{lem:restriction-zero}.

There is also an exact sequence
   \[
   0 \arr \ch{\torusgl}\stackrel{\sigma_{1}}{\longrightarrow}\ch{\torusgl}
      \arr \ch{\torussl} \arr 0,
   \]
so $\ch{\torussl}$ is the quotient $\ch{\torusgl}/(\sigma_{1})$.

\begin{lemma}\label{lem:permutation-rep}
If $G$ is a subgroup of\/ $\s$, the projection $\bigl(\ch{\torusgl}\bigr)^{G} \arr \bigl(\ch{\torussl}\bigr)^{G}$ induces an isomorphism
   \[
   \bigl(\ch{\torussl}\bigr)^{G}/(\sigma_{1}) \simeq
   \bigl(\ch{\torusgl}\bigr)^{G}.
   \]
\end{lemma}

\begin{proof}
This is equivalent to saying that the exact sequence above stays exact after taking $G$-invariants; but we have that $\rH^{1}\bigl(G, \ch{\torusgl}\bigr) = 0$, because $\ch{\torusgl}$ is a torsion-free permutation module under $G$.
\end{proof}

Part~\refpart{prop:describe-cycpgl}{1} comes from the surjectivity of the restriction homomorphism $\bigl(\ch{\torusgl}\bigr)^{\cyc} \arr \bigl(\ch{\torussl}\bigr)^{\cyc}$ and Proposition~\refall{prop:describe-cycgl}{4}.

We construct the splitting $\bigl(\ch{\torussl}\bigr)^{\cyc} \arr \ch{\cycsl}$ by taking the splitting $\bigl(\ch{\torusgl}\bigr)^{\cyc} \arr \ch{\cycgl}$ constructed in the previous section, tensoring it with $\bigl(\ch{\torusgl}\bigr)^{\cyc}/(\sigma_{1})$ over $\bigl(\ch{\torusgl}\bigr)^{\cyc}$, to get a ring homomorphism
   \[
   \bigl(\ch{\torusgl}\bigr)^{\cyc}/(\sigma_{1}) \arr 
      \ch{\cycgl}/(\sigma_{1})
   \]
and using the isomorphisms
   \[
   \bigl(\ch{\torussl}\bigr)^{\cyc} \simeq
   \bigl(\ch{\torusgl}\bigr)^{\cyc}/(\sigma_{1})
   \]
and
   \[
   \ch{\cycgl}/(\sigma_{1}) \simeq \ch{\cycsl}
   \]
constructed above. This proves part~\refpart{prop:describe-cycpgl}{2}. Part~\refpart{prop:describe-cycpgl}{3} follows from Proposition~\refall{prop:describe-cycgl}{2}.

To prove part~\refpart{prop:describe-cycpgl}{5} consider the diagram of restriction homomorphisms
   \[
   \xymatrix{
   \ch{\cycgl} \ar[r]\ar[d]
   & \ch{\cycsl}\ar[r]\ar[d] 
   & \ch{\cycmu} \ar[d]\\
   \bigl(\ch{\torusgl}\bigr)^{\cyc} \ar[r]
   & \bigl(\ch{\torussl}\bigr)^{\cyc} \ar[r]
   & \ch{\mmu} \hsmash{.}
   }
   \]
The surjectivity of the map in the statement follows from Proposition~\refall{prop:describe-cycgl}{5} and from the fact that the first arrow in the bottom row is surjective.

To prove injectivity take an element $u$ of $\ch{\cycsl}$ that maps to $0$ in $\ch{\cycmu}$ and in $\bigl(\ch{\torussl}\bigr)^{\cyc}$. Let $v$ be an element of $\ch{\cycgl}$ mapping to $u$. Since the kernel of the homomorphism $\bigl(\ch{\torusgl}\bigr)^{\cyc} \arr \bigl(\ch{\torussl}\bigr)^{\cyc}$ is generated by $\sigma_{1}$, we can write the image of $v$ in $\bigl(\ch{\torusgl}\bigr)^{\cyc}$ as $\sigma_{1}w$ for some $w \in \bigl(\ch{\torusgl}\bigr)^{\cyc}$. Then the element $v - \phi(\sigma_{1}w)$ maps to $0$ in $\ch{\mmu}$ and in $\bigl(\ch{\torusgl}\bigr)^{\cyc}$; hence, by Proposition~\refall{prop:describe-cycgl}{3}, it is $0$. So $v = \phi(\sigma_{1})\phi(w)$ maps to $0$ in $\bigl(\ch{\torussl}\bigr)^{\cyc}$, as claimed.

Let us prove part~\refpart{prop:describe-cycpgl}{4}. The statement on Chow rings is an immediate consequence of part~\refpart{prop:describe-cycpgl}{5}.

For the cohomology, we will argue as follows. We have a long exact sequence
   \[
   \xymatrix@C+20pt{
   &\cdots\ar[r]& \H[i-1]{\cycsl} \ar[dll] _{\partial} \\
   \H[i-2]{\cycgl} \ar[r]_{\chern_{1}(V)} 
   & \H[i]{\cycgl} \ar[r]
   &\H[i]{\cycsl} \ar[dll]_{\partial}\\
   \H[i-1]{\cycgl} \ar[r] & \cdots
   }
   \]

By Proposition~\ref{prop:chcycgl-hcycgl}, the cycle homomorphism $\ch{\cycgl} \arr \H{\cycgl}$ is an isomorphism. Hence, for each $i$ we have a commutative diagram with exact rows
   \[
   \xymatrix{
   \ch[i-1]{\cycgl} \ar[r]\ar[d] &\ch[i]{\cycgl}\ar[r]\ar[d]
   &\ch[i]{\cycsl}  \ar[r]\ar[d] & 0\\
   \H[2i-2]{\cycgl} \ar[r] &\H[2i]{\cycgl}\ar[r]
   &\H[2i]{\cycsl}  \ar[r] &\H[2i-1]{\cycgl}= 0
   }
   \]
in which the first two columns are isomorphisms. This implies that the third column is also an isomorphism: so the cycle homomorphism $\ch{\cycsl}\arr \H[\even]{\cycsl}$ is an isomorphism. Therefore the homomorphism $\H[\even]{\cycsl} \arr \H[\even]{\torussl} \times \H[\even]{\cycmu}$ is injective.

When $i$ is odd, we have an exact sequence
   \[
   0 = \H[i]{\cycgl} \arr \H[i]{\cycsl}
   \stackrel{\partial}{\longrightarrow} \H[i]{\cycgl}
   \xarr{\chern_{1}(V)} \H[i+2]{\cycgl};
   \]
hence the boundary homomorphism $\partial\colon \H[\odd]{\cycsl} \H[\even]{\cycgl}$ yields an isomorphism of $\H[\odd]{\cycsl}$ with the annihilator of the element $\chern_{1}(V)$ of $\H[\even]{\cycgl} = \ch{\cycgl}$. From the description of the ring $\ch{\cycsl}$ in \refpart{prop:describe-cycgl}{2}, it is easy to conclude that this annihilator is the ideal generated by $\xi$.

Consider a free action of $\cycgl$ on an open subscheme $U$ of a representation. The diagram of embeddings
   \[
   \xymatrix{
   \cycmu \ar[r]\ar[d] &\cyc\times \gm \ar[d]\\
   \cycsl \ar[r] &\cycgl
   }
   \]
induces a cartesian diagram
   \[
   \xymatrix{
   U/\cycmu \ar[r]\ar[d] &U/\cyc\ltimes \gm \ar[d]\\
   U/\cycsl \ar[r] &U/\cycgl
   }
   \]
in which the rows are principal $\gm$-bundles, and the columns are $\gm$-equivariants. This in turn induces a commutative diagram
   \[
   \xymatrix{
   \H[\odd]{\cycsl} \ar[r]^{\partial}\ar[d]
   & \H[\even]{\cycgl} \ar[d] \\
   \H[\odd]{\cycmu}\ar[r]
   & \H[\even]{\cyc \times \gm}
   }
   \]
in which the top row is injective, and has as its image the ideal $(\xi) \subseteq \H[\even]{\cycgl}$ as we have just seen. Furthermore, every element of $(\xi) \subseteq \H[\even]{\cycgl}$ maps to $0$ in $\H{\torusgl}$, because it is torsion: hence $(\xi)$ injects into $\H[\even]{\cycmu}$, by Proposition \refall{prop:describe-cycgl}{3}. Since $\cycmu$ is contained into $\cyc\times\gm$, it follows that $(\xi)$ also injects into $\H[\even]{\cyc\times\gm}$. So the composite arrow $\H[\odd]{\cycsl} \arr \H[\even]{\cyc\times\gm}$ in the commutative diagram above is injective. It follows that the left hand column is injective.

This ends the proof of Proposition~\ref{prop:describe-cycpgl}.
\end{proof}

\begin{proof}[Proof of Proposition~\ref{prop:ker-whole}]
We need to analyze the action of the normalizer $\N$ of $\cyc$ in $\s$ on the Chow ring $\ch{\cycpgl}$. If we identify $\{1, \dots, p\}$ with the field $\mathbb{F}_p$ with $p$ elements, by sending each $i$ into its class modulo $p$, then $\cyc$ can be identified with the additive group $\mathbb{F}_p$ itself, acting by translations. There is also the multiplicative subgroup $\mathbb{F}^*_p$ of $\s$, acting via multiplication. This is contained in the normalizer of $\cyc = \mathbb{F}_p$, and, since $p$ is a prime, it is easy to show that the normalizer of $\cyc$ inside $\s$ is in fact the subgroup generated by $\mathbb{F}_p$ and $\mathbb{F}^*_p$, which is the semi-direct product $\norm$.

The subgroup $\cyc = \FF_{p}$ acts trivially, so in fact the action is through $\FF_{p}^{*}$. The action of $\FF_{p}^{*}$ leaves $\mmu$ invariant, and the result of the action of $a \in \FF_{p}^{*}$ on $\zeta \in \mmu$ is $\zeta^{a}$: hence $a$ acts on $\ch{\mmu} = \ZZ[\eta]/(p\eta)$ by sending $\eta$ to $a\eta$, and the ring of invariants is the subring generated by $\eta^{p-1}$. The image of $\ch{\toruspgl}$ into $\ch{\mmu}$ is the subring generated by $\eta^{p}$, by Proposition~\ref{prop:describe-cycpgl}, and its intersection with the ring of invariants in $\ch{\mmu}$ is the subring generated by $\eta^{p(p-1)}$. This shows that the image of $\bigl(\ch{\toruspgl}\bigr)^{\s}$ into $\ch{\mmu}$ is contained in the subring generated by $\eta^{p(p-1)}$. The opposite inclusion is ensured by the fact that the discriminant $\delta \in \bigl(\ch{\toruspgl}\bigr)^{\s}$ maps to $-\eta^{p(p-1)}$.
\end{proof}

\section{On $\sympgl$}

The group $\s$ does not act on $\cycpgl$, only the normalizer $\norm$ of $\cyc$ does. Nevertheless, we define the subring $\bigl(\ch{\cycpgl}\bigr)^{\s}$ of $\ch{\cycpgl}$ consisting of all the elements that are  invariant under $\norm$, and whose images in $\ch{\toruspgl}$ are $\s$-invariant. The restriction homomorphism $\ch{\sympgl} \to \ch{\cycpgl}$ has its image in $\bigl(\ch{\cycpgl}\bigr)^{\s}$.

The result we need about  $\sympgl$ is the following. 

\begin{proposition}\label{prop:on-sympgl}
The localized restriction homomorphism
   \[
   \ch{\sympgl}\invert \longrightarrow \bigl(\ch{\cycpgl}\bigr)^{\s}\invert
   \]
is an isomorphism.
\end{proposition}

Of course the statement can not be correct without inverting $(p-1)!$, because the torsion part of $\ch{\cycpgl}$ is all $p$-torsion, while $\ch{\sympgl}$ contains a lot of $(p-1)!$-torsion coming from $\ch{\s}$. This is complicated, but fortunately we do not need to worry about it.

\begin{proof}
Injectivity is clear: because of the projection formula, the composite
   \[
   \tsf{\sympgl}{\cycpgl} \res{\cycpgl}{\sympgl}
   \colon \ch{\sympgl} \arr \ch{\cycpgl}
   \]
is multiplication by $\tsf{\sympgl}{\cycpgl} = (p-1)!$.

To show surjectivity, take a class $u \in \bigl(\ch{\cycpgl}\bigr)^{\s}$, and set 
   \[
   v \eqdef \tsf{\cycpgl}{\sympgl} u \in \ch{\sympgl}.
   \]
We apply  Mackey's formula (Proposition~\ref{prop:mackey}). The double quotient
   \[
   {\cycpgl}\backslash{\sympgl}/{\cycpgl} = {\cyc}\backslash{\s}/{\cyc}
   \]
consists of $p-1$ elements coming from the normalizer $\norm$, and $(p-1) \frac{(p-2)! - 1}{p}$ elements with the property that, if we call $s$ a representative in $\sympgl$, we have
   \[
   s (\cycpgl) s ^{-1} \cap \cycpgl = \toruspgl.
   \]
Therefore
   \[
   \res{\sympgl}{\cycpgl} v = (p-1)u +
   (p-1)\frac{(p-2)! - 1}{p} \tsf{\toruspgl}{\cycpgl} 
   \res{\cycpgl}{\toruspgl}u;
   \]
hence it is enough to show that an element in the image of the trasfer map 
   \[
   \tsf{\toruspgl}{\cycpgl} \colon \invtorus
   \longrightarrow \bigl(\ch{\cycpgl}\bigr)^{\s}
   \]
is in the image of $\ch{\toruspgl}$, up to a multiple of $(p-1)!$. But again an easy application of Mackey's formula reveals that
   \[
    \res{\cycpgl}{\toruspgl}\tsf{\toruspgl}{\sympgl} w = (p-1)!w
   \]
for all $w \in \invtorus$, and this finishes the proof.
\end{proof}

\section{Some results on $\ch{\pgl}$}\label{sec:aux-pgl}

In this section we prove some auxilliary results, which play an important role in the proof of the main theorems.

The following observation is in \cite[Corollary~2.4]{vezzosi-pgl3}.

\begin{proposition}\label{prop:p-torsion}
If $\xi$ is a torsion element of $\ch{\pgl}$, or $\H{\pgl}$, then $p\xi = 0$.
\end{proposition}

\begin{proof}
Suppose that $\xi \in \ch[m]{\pgl}$. Take a representation $V$ of $\pgl$ with an open subset $U$ on which $\pgl$ acts freely, such that the codimension of $V \setminus U$ has codimension larger than $m$, so that $\ch[m]{\pgl} = \ch[m]{}(B)$, where we have set $B \eqdef U/\pgl$. Let $\pi\colon E \arr B$ be the Brauer--Severi scheme associated with the $\pgl$-torsor $U \arr B$: this is the projection $U/H \arr U/\pgl$, where $H$ is the parabolic subgroup of $\pgl$ consisting of classes of matrices $(a_{ij})$ with $a_{i1} = 0$ when $i > 1$.
The embedding $H \into \pgl$ lifts to an embedding $H \into \GL$, as the subgroup of matrices $(a_{ij})$ with $a_{i1} = 0$ when $i > 1$, and $a_{11} = 1$; hence the pullback $\ch[m]{}(B) \arr \ch[m]{}(E)$ factors through $\ch[m]{\GL}$, which is torsion-free. It follows that $\xi$ maps to $0$ in $\ch[m]{}(E)$.

Now consider the Chern class $\chern_{p-1}(\rT_{E/B}) \in \ch[p-1]{}(E)$ of the relative tangent bundle. This has the property that $\pi_{*}\chern_{p-1}(\rT_{E/B}) = p[B] \in \ch[0]{}(B)$; hence, by the projection formula we have
   \begin{align*}
   p\xi &= \xi\cdot\pi_{*}\bigl(\chern_{p-1}(\rT_{E/B})\bigr) \\
   &= \pi_{*}\bigl(\pi^{*}\xi \cdot\chern_{p-1}(\rT_{E/B})\bigr)\\
   &= 0.
   \end{align*}

The proof for cohomology is identical, except for notation.
\end{proof}

\begin{proposition}
The restriction homomorphisms $\ch{\pgl} \arr \ch{\cycpgl}$ and $\H{\pgl} \arr \H{\cycpgl}$ are injective.
\end{proposition}

\begin{proof}
By a classical result of Gottlieb (\cite{gottlieb-euler}) the homomorphism $\H{\pgl} \arr \H{\sympgl}$ is injective; while the injectivity of $\ch{\pgl} \arr \ch{\sympgl}$ is a recent result of Totaro. This is unpublished: a sketch of proof is presented in \cite{vezzosi-pgl3}.

\begin{theorem}[Totaro]\label{thm:totaro-injectivity}
If $G$ is a connected linear algebraic group over a field $k$ of characteristic $0$ acting on a scheme $X$ of finite type over $k$, and $N$ is the normalizer of maximal torus, then the restriction homomorphism $\ch{G}(X) \arr \ch{N}(X)$ is injective.
\end{theorem}

Now, the kernels of the homomorphisms in the statement are $p$-torsion, by Proposition~\ref{prop:p-torsion}, while the kernels of
   \[
   \ch{\sympgl} \arr \ch{\cycpgl}\quad
   \text{and} \quad\H{\sympgl} \arr \H{\cycpgl}
   \]
are $(p-1)!$-torsion, by the projection formula, so the statement follows.
\end{proof}

Here is the basic result that we are going to use in order to verify that a given relation holds in $\ch{\pgl}$ and $\H{\pgl}$.

\begin{proposition}\label{prop:keyinjectivity}
The homomorphisms
   \[
   \ch{\pgl} \arr \ch{\toruspgl} \times \ch{\cycmu}
   \]
and
   \[
   \H{\pgl} \arr \H{\toruspgl} \times \H{\cycmu}
   \]
obtained from the embeddings $\toruspgl \into \pgl$ and $\cycmu \into \pgl$ are injective.
\end{proposition}

\begin{proof}
This follows from Propositions \ref{prop:keyinjectivity} and \refall{prop:describe-cycpgl}{4}.
\end{proof}

Here is another fundamental fact, which is one of the cornerstones of the treatment of $\pgl[3]$ in \cite{vezzosi-pgl3}. In the Lie algebra $\liesl$ of matrices of trace $0$ consider the Zariski open subset $\liesls$ consisting of matrices with distinct eigenvalues; this is invariant by the action of $\pgl$. Furthermore, we will consider the subspace $\diag \subseteq \liesl$ of diagonal matrices with trace equal to zero, and $\diags = \diag \cap \liesls$. The subspaces $\diag$ and $\diags$ are invariant under the action of $\sympgl \subseteq \pgl$.

\begin{proposition}[see \cite{vezzosi-pgl3}, Proposition~3.1] \label{prop:rest-isom}
The composites of restriction homomorphisms
   \[
   \ch{\pgl}(\liesls) \longrightarrow \ch{\sympgl}(\liesls)
   \longrightarrow \ch{\sympgl}(\diags)
   \]
and
   \[
   \H{\pgl}(\liesls) \longrightarrow \H{\sympgl}(\liesls)
   \longrightarrow \H{\sympgl}(\diags)
   \]
are isomorphisms.
\end{proposition}

\begin{proof}
The $\sympgl$-equivariant embedding $\diags \subseteq \liesls$ induces a $\pgl$\dash equivariant morphism $\pgl \times^{\sympgl} \diags \to \liesls$, which sends the class of a pair $(A, X)$ into $AXA^{-1}$. This morphism is easily seen to be an isomorphism, and the proof follows.
\end{proof}

\begin{corollary}\label{cor:mapstorestriction-equal}
The restriction homomorphisms
   \[
   \ch{\pgl} \to \ch{\toruspgl} \quad\text{and} \quad
   \ch{\sympgl} \to \ch{\toruspgl}
   \]
have the same image.
\end{corollary}

\begin{proof}
In the commutative diagram of restriction homomorphisms
   \[
   \xymatrix{
   \ch{\pgl} \ar[r]\ar[d] & \ch{\pgl}(\liesls)\ar[d]\\
   \ch{\toruspgl} \ar[r]  & \ch{\toruspgl}(\liesls)
   }
   \]
the top row is surjective. On the other hand, the action on $\toruspgl$ on $\liesls$ is trivial and $\liesls$ is an open subscheme of an affine space, so the bottom row is an isomorphism. It follows that the image of $\ch{\pgl}$ in $\ch{\toruspgl}$ maps isomorphically onto the image of $\ch{\pgl}(\liesls)$ in $\ch{\toruspgl}(\liesls)$. A similar argument shows that the image of $\ch{\sympgl}$ in $\ch{\toruspgl}$ maps isomorphically onto the image of $\ch{\sympgl}(\diags)$ in $\ch{\toruspgl}(\diags)$. By we also have a commutative diagram
   \[
   \xymatrix{
   \ch{\pgl}(\liesls) \ar[r]\ar[d] & \ch{\sympgl}(\diags)\ar[d]\\
   \ch{\toruspgl}(\liesls) \ar[r]  & \ch{\toruspgl}(\diags)
   }
   \]
where the top row is an isomorphism, by Proposition~\ref{prop:rest-isom}, and this concludes the proof.
\end{proof}

\section{Localization}

Consider the top Chern classes
   \[
   \chern_{p^2-1}(\liesl) \in \ch[p^2-1]{\pgl}
   \quad\text{and}\quad
   \chern_{p-1}(\diag) \in \ch[p-1]{\sympgl}.
   \]We have the following fact.

\begin{proposition}\label{prop:localization}
The restriction homomorphism $\ch{\pgl} \to \ch{\sympgl}$ carries $\chern_{p^2-1}(\liesl)$ into the ideal $\bigl(\chern_{p-1}(\diag)\bigr) \subseteq\ch{\sympgl}$. The induced homomorphism
   \[
   \ch{\pgl}\big/\bigl(\chern_{p^2-1}(\liesl)\bigr) \longrightarrow
   \ch{\sympgl}/\bigl(\chern_{p-1}(\diag)\bigr)
   \]
becomes an isomorphism when tensored with $\mathbb{Z}[1/(p-1)!]$.
\end{proposition}

\begin{proof}
The representation $\diag$ of $\sympgl$ is naturally embedded in $\liesl$,
so we have that
   \[
   \chern_{p^2-1}(\liesl) = \chern_{p-1}(\diag)
   \chern_{p+1}(\liesl/\diag) \in \ch[p^{2}-1]{\sympgl},
   \]
and this proves the first statement.

The pullbacks
   \[
   \ch{\pgl} \arr \ch{\pgl}(\liesl \setminus \{0\}) \quad \text{and}\quad 
   \ch{\sympgl} \arr \ch{\sympgl}(\diag \setminus \{0\})
   \]
are surjective, and their kernels are the ideals generated by $\chern_{p^2-1}(\liesl)$ and $\chern_{p-1}(\diag)$ respectively: so it enough to show that the homomorphism
   \[
   \ch{\pgl}(\liesl \setminus \{0\}) \arr \ch{\sympgl}(\diag \setminus \{0\})
   \]
obtained by restricting the groups, and then pulling back along the embedding $\diag \setminus \{0\} \into \liesl \setminus \{0\}$ becomes an isomorphism after inverting $(p-1)!$.

Now, consider the diagram
   \[
   \xymatrix{
   {}\ch{\pgl}(\liesl \setminus \{0\}) \ar[r] \ar[d]&
      {}\ch{\pgl}(\liesl^{*}) \ar[d]\\
   {}\ch{\sympgl}(\diag \setminus \{0\}) \ar[r]&
      {}\ch{\sympgl}(\diag^{*})
   }
   \]
where all the arrows are the obvious ones. The rows are surjective, while the right hand column is an isomorphism, by Proposition~\ref{prop:rest-isom}: hence it is enough to show that the rows are injective, after inverting $(p-1)!$.

The first step is to observe that the restriction homomorphism $\ch{\pgl}(\liesl \setminus \{0\}) \arr \ch{\sympgl}(\liesl \setminus \{0\})$ is injective, by Totaro's Theorem~\ref{thm:totaro-injectivity}. Next, the restriction homomorphisms $\ch{\sympgl}(\liesl \setminus \{0\}) \arr \ch{\cycpgl}(\liesl \setminus \{0\})$ and $\ch{\sympgl}(\diag^{*}) \arr \ch{\cycpgl}(\diag^{*})$ become injective after inverting $(p-1)!$. So it is enough to show that the restriction homomorphisms $\ch{\cycpgl}(\liesl \setminus \{0\}) \arr \ch{\cycpgl}(\liesl^{*})$ and $\ch{\cycpgl}(\diag \setminus \{0\}) \arr \ch{\cycpgl}(\diag ^{*})$ are injective.

\begin{lemma}\label{lem:localization}
Suppose that $W$ is a representation of $\cycpgl$, and $U$ an open subset of $W\setminus \{0\}$. Assume that

\begin{enumeratea}

\item the restriction of $W$ to $\cycmu$ splits as a direct sum of $1$-dimensional representations $W = L_{1} \oplus \dots \oplus L_{r}$, in such a way that the characters $\cycmu \arr \gm$ describing the action of $\cycmu$ on the $L_{i}$ are all distinct, and each $L_{i}\setminus \{0\}$ is contained in $U$, and

\item $U$ contains a point that is fixed under $\toruspgl$.

\end{enumeratea}

Then the restriction homomorphism $\ch{\cycpgl}(W \setminus \{0\}) \arr \ch{\cycpgl}(U)$ is an isomorphism.
\end{lemma}

\begin{proof}
First of all, let us show that $\ch{\cycmu}(W \setminus \{0\}) \arr \ch{\cycmu}(U)$ is an isomorphism. Denote by $D$ the complement of $U$ in $W \setminus \{0\}$, with its reduced scheme structure. Let $P$ be the projectivization of $W$, and call $\overline{U}$ and $\overline{D}$ the (respectively open and closed) subschemes of $P$ corresponding to $U$ and $D$. We have a commutative diagram
   \[
   \xymatrix{
   {}\ch{\cycmu}(\overline{D}) \ar[r]\ar[d]
      &{}\ch{\cycmu}(P) \ar[r]\ar[d]
      &{}\ch{\cycmu}(\overline{U}) \ar[r] \ar[d]
      & 0\\
   {}\ch{\cycmu}(D) \ar[r]
      &{}\ch{\cycmu}(W \setminus \{0\}) \ar[r]
      &{}\ch{\cycmu}(U) \ar[r]
      & 0
   }
   \]
where the columns are surjective pullbacks, and the rows are exact. It follows that is enough to show that the composite
   \[
   \ch{\cycmu}(\overline{D}) \arr \ch{\cycmu}(P) 
   \arr \ch{\cycmu}(W \setminus \{0\})
   \]
is $0$, or, equivalently, that any element of the kernel of $\ch{\cycmu}(P) \arr \ch{\cycmu}(\overline{U})$ maps to $0$ in $\ch{\cycmu}(W \setminus \{0\})$. Denote by $q_{i} \in P$ the rational point corresponding to $L_{i}$.

Denote by $\ell_{i} \in \ch[1]{\cycmu}$ the first Chern class of the character $\cycmu \arr \gm$ describing the action of $\cycmu$ on $L_{i}$, and $h \in \ch[1]{\cycmu}$ the first Chern class of the sheaf $\cO(1)$ on $P$. We have presentations
   \[
   \ch{\cycmu}(P) =  \ZZ[\xi, \eta, h]/
   \bigl(p\xi, p\eta, (h-\ell_{1}) \dots (h-\ell_{r})\bigr)
   \]
and
   \[
   \ch{\cycmu}(W \setminus \{0\}) = \ZZ[\xi, \eta]/
   (p\xi, p\eta, \ell_{1} \dots \ell_{r}),
   \]
and a commutative diagram
   \[
   \xymatrix{
   {}\ZZ[\xi, \eta, h]/ \bigl(p\xi, p\eta, (h-\ell_{1}) \dots
      (h-\ell_{r})\bigr) \ar[r]\ar[d]& 
   {}\ZZ[\xi, \eta]/ (p\xi, p\eta, \ell_{1} \dots \ell_{r}) \ar[d]\\
   {}\FF_{p}[\xi, \eta, h]/ \bigl((h-\ell_{1}) \dots (h-\ell_{r})\bigr) \ar[r]
   & \FF_{p}[\xi, \eta]/ (\ell_{1} \dots \ell_{r})
   }
   \]
in which the first row is the map that sends $h$ to $0$, and corresponds to the pullback.

The restriction homomorphism $\ch{\cycmu}(P) \arr \ch{\cycmu}(q_{i}) = \ch{\cycmu}$ sends $h$ into $\ell_{i}$. But $\overline{U}$ contains all the $q_{i}$, so the kernel $K$ of the restriction $\ch{\cycmu}(P) \arr \ch{\cycmu}(\overline{U})$ is contained in the intersection of the ideals $(h - \ell_{i})$. In the polynomial ring $\FF_{p}[\xi, \eta, h]$, however, the intersection of the ideals $(h - \ell_{i})$ is the ideal generated by the product of the $h - \ell_{i}$, because $\FF_{p}[\xi, \eta, h]$ is a unique factorization domain, and the $h - \ell_{i}$ are pairwise non-associated primes. Hence the image of an element of $K$ is $0$ in $\FF_{p}[\xi, \eta]/ (\ell_{1} \dots \ell_{r})$; but the homomorphism
   \[
   \ch{\cycmu}(W \setminus \{0\})
   \arr \FF_{p}[\xi, \eta]/ (\ell_{1} \dots \ell_{r})
   \]
   is an isomorphism in positive degree, and from this the statement follows.

Now consider the restriction homomorphism
   \[
   \ch{\cycpgl}(W \setminus \{0\}) \arr \ch{\cycpgl}(U).
   \]
Denote by $\gamma$ the top Chern class of $W$ in $\ch{\cycpgl}$; the kernel of the surjective pullback $\ch{\cycpgl} \arr \ch{\cycpgl}(W \setminus \{0\})$ is the ideal generated by $\gamma$, and we need to show that the kernel of the pullback $\ch{\cycpgl} \arr \ch{\cycpgl}(U)$ is also the ideal generated by $\gamma$.

Denote by $R$ the image of $\ch{\cycpgl}$ in $\ch{\cycmu} = \ZZ[\xi, \eta]/(p\xi, p\eta)$; this is the subring generated by $\xi$ and the image of $\sigma_{p}$, that is $\eta^{p} - \xi^{p-1}\eta$, by Lemma~\ref{lem:restrictcanonical-mmu}.

Take some $u$ in the kernel of $\ch{\cycpgl} \arr \ch{\cycpgl}(U)$. Since $\toruspgl$ has a fixed point in $U$, the pullback $\ch{\toruspgl} \arr \ch{\toruspgl}(U)$ is an isomorphism; hence $u$ is contained in the kernel of the restriction $\ch{\cycpgl} \arr \ch{\toruspgl}$. This kernel is the ideal $\xi\ch{\cycpgl}$, which is a vector space over the field $\FF_{p}$, with a basis consisting of the elements $\xi^{i}\sigma_{p}^{j}$, with $i > 0$ and $j \geq 0$. The homomorphism $\ch{\cycpgl} \arr \ch{\cycmu}$ sends $\xi^{i}\sigma_{p}^{j}$ into $\xi^{i}(\eta^{p} - \xi^{p-1}\eta)^{j}$. The two elements $\xi$ and $\eta^{p} - \xi^{p-1}\eta$ are algebraically independent in $\ch{\cycmu}$, so the ideal $\xi\ch{\cycpgl}$ map isomorphically onto the ideal $\xi R$. Hence it is enough to show that $u$ maps into the ideal $\gamma R$. But $u$ maps into the ideal $\gamma \ch{\cycmu}$, because by hypothesis maps into $0$ in $\ch{\cycmu}(W \setminus \{0\})$, so we will be done once we have shown that $\gamma R = R \cap \gamma \ch{\cycmu}$.

For this purpose, consider the diagram
   \[
   \xymatrix{
   R \ar@{ >->}[r] \ar[d]& {}\ch{\cycmu}\ar[d]\\
   {}\FF_{p}[\xi, \eta^{p} - \xi^{p-1}\eta] \ar@{ >->}[r]
      &\FF_{p}[\xi, \eta]
   }
   \]
where the horizontal arrows are inclusions and the vertical arrows are isomorphisms in positive degree. It suffices to prove that
   \[
   \gamma \FF_{p}[\xi, \eta^{p} - \xi^{p-1}\eta] = 
      \FF_{p}[\xi, \eta^{p} - \xi^{p-1}\eta] \cap \gamma\FF_{p}[\xi, \eta];
   \]
but this follows from the fact that the extension $\FF_{p}[\xi, \eta^{p} - \xi^{p-1}\eta] \subseteq \FF_{p}[\xi, \eta]$ is faithfully flat, since it is a finite extension of regular rings.

This concludes the proof of the Lemma.
\end{proof}
The Lemma applies to the case $W = \diag$ and $W = \liesl$. In the first case this is straightforward; in the second one it follows from Lemma~\ref{lem:describe-pgl}.
\end{proof}

\section{The classes $\rho$ and $\beta$}\label{sec:rho}

In this section we construct the classes $\rho \in \ch[p+1]{\pgl}$ and $\beta \in \H[3]{\pgl}$. 

\begin{proposition}
There exists a unique torsion class $\rho \in \ch[p+1]{\pgl}$, whose image in $\ch[p+1]{\cycmu}$ equals $r = \xi\eta(\xi^{p-1} - \eta^{p-1})$.

Furthermore we have $\rho^{p-1} = \chern_{p^{2}-1}(\liesl) \in \ch{\pgl}$.
\end{proposition}

\begin{proof}
Uniqueness is obvious from Proposition~\ref{prop:keyinjectivity}.

Let us construct a $p$-torsion element $\overline{\rho} \in \ch[p+1]{\sympgl}$ that maps to $r$ in $\ch{\cycmu}$.

Consider the element $-\xi\phi(\sigma_{p}) = \xi\chern_{p}(V) \in \ch[p+1]{\cycpgl}$; by  Lemma~\ref{lem:restrictcanonical-mmu}, its restriction to $\ch{\cycmu}$ is $r$. It is $p$-torsion, because $\xi$ is $p$-torsion; hence it maps to $0$ in $\ch{\toruspgl}$. Since the torsion part of $\ch{\cycmu}$ injects into $\ch{\cycmu}$, and the image of $\xi\chern_{p}(V)$ in $\ch{\cycmu}$ is invariant under $\norm$, it follows that $\xi\chern_{p}(V)$ is also invariant under $\norm$.

By Proposition~\ref{prop:on-sympgl}, there exists a $p$-torsion class $\overline{\rho} \in \ch[p+1]{\sympgl}$ whose image in $\ch{\cycmu}$ is $r$. By Proposition~\ref{prop:localization}, there exists a $p$-torsion element $\rho \in \ch[p+1]{\pgl}$ whose image in $\ch[p+1]{\sympgl}$ has the form $\overline{\rho} + \chern_{p-1}(D_{p})\sigma$ for a certain class $\sigma \in \ch[2]{\sympgl}$.

The image of $\rho$ in $\bigl(\ch{\cycmu}\bigr)^{\slfinite} = \ZZ[q,r]$ must be an integer multiple $ar$ of $r$, for reasons of degree. The image of $\chern_{p-1}(D_{p})$ is $-\xi^{p-1}$; hence by mapping into $\ch{\cycmu}$ we get an equality
   \[
   ar = r - \xi^{p-1}h \in \ch{\cycmu},
   \]
where $h \in \ch[2]{\cycmu}$ is the image of $\sigma$. From this equality it follows easily that $a$ is $1$ and $h$ is $0$, and therefore $\rho$ maps to $r$.

To check that $\rho^{p-1} = \chern_{p^{2}-1}(\liesl)$, observe that both members of the equality are 0 when restricted to $\toruspgl$; hence, by Proposition~\ref{prop:keyinjectivity}, it is enough to show that the restriction of $\chern_{p^{2}-1}(\liesl) = \chern_{p^{2}-1}(\liegl)$ to $\ch[p^{2}-1]{\cycmu}$ equals $r^{p-1}$; and this follows from Lemma~\ref{lem:restrict-chern}.
\end{proof}

\begin{corollary}
The restriction homomorphism $\ch{\pgl} \arr \bigl(\ch{\cycmu}\bigr)^{\slfinite}$ is surjective.
\end{corollary}

\begin{proof}
The ring $\bigl(\ch{\cycmu}\bigr)^{\slfinite}$ is generated by $q$ and $r$. The class $-\chern_{p^{2}-p}(\liesl)$ restricts to $q$, by Lemma~\ref{lem:restrict-chern}, while $\rho$ restricts to $r$.
\end{proof}

\begin{remark}\label{rmk:azumaya}
The class $\rho$ gives a new invariant for sheaves of Azumaya algebras of prime degree. Let $X$ be a scheme of finite type over $k$, and let $\cA$ be a sheaf of  Azumaya algebras of degree $p$. This corresponds to a $\pgl$-torsor $E \arr X$; and according to a result of Totaro (see \cite{totaro-classifying} and \cite{edidin-graham-characteristic}), we can associate to the class $\rho \in \ch[p+1]{\pgl}$ and the $\pgl$-torsor $E$ a class $\phi(\cA) \in \ch[p+1]{}(X)$ (where by $\ch{}(X)$ we mean the bivariant ring of $X$, see \cite{fulton}). Since by definition $\cA$ is the vector bundle associated with $E$ and the representation $\liegl$ of $\pgl$, we have the relation
   \[
   \rho(\cA)^{p-1} = \chern_{p^{2}-1}(\cA) \in \ch[p^{2}-1]{}(X).
   \]
\end{remark}

\begin{remark}
The class $\rho$ depends on the choice of the primitive $p\th$ root of $1$ that we have denoted by $\omega$. If we substitute $\omega^{i}$ for $\omega$, then the new class $\rho$ is $i\rho$.
\end{remark}

For the class $\beta$, one possibility it to obtain it as the Brauer class of the canonical $\pgl$-principal bundle, as explained in the Introduction. Another possibility is to define it via a transgression homomorphism, as follows. There is a well known Hochshild--Serre spectral sequence
   \[
   E_{2}^{ij} = \H[i]{\pgl}\otimes\H[j]{\gm} \Longrightarrow
      \H[i+j]{\GL}
   \]
from which we get an exact sequence
   \[
   \H[2]{\GL} \arr \H[2]{\gm} \arr \H[3]{\pgl} \arr \H[3]{\toruspgl} = 0;
   \]
and $\H[2]{\gm}$ is the infinite cyclic group generated by the first Chern class $t$ of the identity character $\gm = \gm$, while $\H[2]{\GL}$ is the cyclic group generated by the first Chern class of the determinant $\GL = \gm$, whose image in $\H[2]{\gm}$ is $pt$. Hence $\H[3]{\pgl}$ is the cyclic group of order $p$ generated by the image of $t$. We define $\beta \in \H[3]{\pgl}$ to be this image.

The odd dimensional cohomology $\H[\odd]{\pgl}$ maps to $0$ in $\H{\toruspgl}$; hence, according to Proposition~\ref{prop:keyinjectivity}, maps injectively into $\H{\cycmu}$. By the results of Section~\ref{sec:cycmu}, we have that $\H[3]{\cycmu}$ is isomorphic to $\ZZ/p\ZZ$, hence the restriction homomorphism $\H[3]{\pgl} \arr \H[3]{\cycmu}$ is an isomorphism; and the image of $\beta$ generated $\H[3]{\cycmu}$. From Proposition~\ref{prop:H-invariants-cycmu} we obtain the following.

\begin{corollary}
The restriction homomorphism $\H{\pgl} \arr \bigl(\H{\cycmu}\bigr)^{\slfinite}$ is surjective.
\end{corollary}

\section{The splitting}\label{sec:splitting}

In this section we prove Theorem~\ref{thm:main-splitting}.

Consider the embeddings
   \[
   \xymatrix{
   \mmu \ar@{^(->}[r]\ar@{^(->}[d]& \toruspgl\ar@{^(->}[d]\\
   \cycmu\ar@{^(->}[r] & \sympgl
   }
   \]
which induce a diagram of restriction homomorphisms
   \[
   \xymatrix{
   \ch{\sympgl} \ar[r]\ar[d] & \invtorus\ar[d]\\
   \ch{\cycmu}\ar[r]     & \ch{\mmu}\hsmash{.}
   }
   \]

\begin{lemma}\label{lem:ontoker-sympgl}
The induced homomorphism
   \[
   \ker\bigl(\ch{\sympgl} \to \ch{\cycmu}\bigr) \longrightarrow
   \ker\bigl(\invtorus \to \ch{\mmu}\bigr)
   \]
is surjective.
\end{lemma}

\begin{proof}
We will prove surjectivity in two steps; first we will show that the map is surjective when tensored with $\mathbb{Z}[1/p]$, then that is surjective when tensored with $\mathbb{Z}[1/(p-1)!]$.

For the first case, notice that $\ch{\cycmu}\invertp$ is 0 in positive degree, while in degree~0 there is nothing to prove; so what we are really trying to show is that $\ch{\sympgl}\invertp \to \invtorus\invertp$ is surjective.

Consider the subgroup $\s[p-1] \subseteq \s$ of the Weyl group of $\pgl$, consisting of permutations of $\{1, \dots, p\}$ leaving $p$ fixed.

\begin{lemma}
The restriction homomorphism $\ch{\s[p-1] \ltimes \toruspgl} \to (\ch{\toruspgl})^{\s[p-1]}$ is surjective.
\end{lemma}

\begin{proof}
There is an isomorphism $\torusgl[p-1] \simeq \toruspgl$, defined by
   \[
   (t_1, \dots, t_{p-1}) \mapsto (t_1, \dots, t_{p-1}, 1)
   \]
that is $\s[p-1]$-equivariant, and therefore induces an isomorphism of the semi-direct product $\s[p-1] \ltimes \toruspgl$ with the normalizer $\symgl[p-1]$ of the maximal torus in $\GL[p-1]$. Hence it is enough to show that $\ch{\symgl[p-1]} \to (\ch{\torusgl[p-1]})^{\s[p-1]}$ is surjective; but the composite
   \[
   \ch{\GL[p-1]} \arr \ch{\symgl[p-1]} \arr (\ch{\torusgl[p-1]})^{\s[p-1]}
   \]
is an isomorphism, and this proves what we want.
\end{proof}

Take an element $u \in \invtorus$; according to the Lemma above, there is some $v \in \ch{\s[p-1] \ltimes \toruspgl}$ such that $\res{\s[p-1] \ltimes \toruspgl}{\toruspgl} v = u$. Consider the element
   \[
   w \eqdef \tsf{\s[p-1] \ltimes \toruspgl}{\sympgl}v \in \ch{\sympgl};
   \]
to compute its restriction to $\ch{\toruspgl}$ we use Mackey's formula (Proposition~\ref{prop:mackey}). The double quotient $\toruspgl \backslash \sympgl / \s[p-1] \ltimes \toruspgl$ has $p$ element, and we may take $\cyc$ as a set of representatives. Then the formula gives us that the restriction of $w$ to $\ch{\toruspgl}$ is
   \[
   \sum_{s \in \cyc} s\res{\s[p-1] \ltimes \toruspgl}{\toruspgl} v = pu.
   \]
If we invert $p$, this shows that $u$ is in the image of $\ch{\sympgl}$, and completes the proof of the first step.

For the second step, take some $u \in L$. According to Proposition \refall{prop:describe-cycpgl}{5} there exists $v$ in the kernel of the restriction homomorphism $\ch{\cycpgl} \to \ch{\cycmu}$ whose restriction to $\ch{\toruspgl}$ is $u$. Consider the element
   \[
   w \eqdef \tsf{\cycpgl}{\sympgl}v.
   \]
I claim that $w$ is in $K$. In fact the restriction of $w$ to $\ch{\toruspgl}$ is $(p-1)!v = -v$, and thefore further restricting it to $\cycpgl$ sends it to $0$.

The double quotient $\toruspgl \backslash \sympgl / \s[p-1] \ltimes \toruspgl$ has $(p-1)!$ elements, and a set of representatives is given by $\s[p-1]$. Hence according to Mackey's formula we have that the restriction of $w$ to $\ch{\toruspgl}$ is
   \[
   \sum_{s \in \s[p-1]} s\res{\cycpgl}{\toruspgl} v = (p-1)!u
   \]
and this completes the second step in the proof of Lemma~\ref{lem:ontoker-sympgl}.
\end{proof}

Similarly, there is a diagram of restriction homomorphisms
   \[
   \xymatrix{
   \ch{\pgl} \ar[r]\ar[d] & \invtorus\ar[d]\\
   \ch{\cycmu}\ar[r]     & \ch{\mmu}\hsmash{.}
   }
   \]

\begin{lemma}\label{lem:ontoker-pgl}
The homomorphism
   \[
   \ker\bigl(\ch{\pgl} \to \ch{\cycmu}\bigr) \longrightarrow
   \ker\bigl(\invtorus \to \ch{\mmu}\bigr)
   \]
induced by restriction is an isomorphism.
\end{lemma}

\begin{proof}
Injectivity follows from Proposition~\ref{prop:keyinjectivity}.

As in the previous case, we show surjectivity first after inverting $p$, and then after inverting $(p-1)!$. 

As before, we have $\ch{\cycmu}\invertp = \ZZ[1/p]$, so we only need to check that $\ch{\pgl}\invertp \arr \invtorus\invertp$ is surjective. This follows from Lemma~\ref{lem:ontoker-sympgl} and from Corollary~\ref{cor:mapstorestriction-equal}.

Now we invert $(p-1)!$. Choose an element
   \[
   u \in \ker\Bigl(\invtorus\to \ch{\mmu}\Bigr)\invert;
   \]
by Lemma~\ref{lem:ontoker-sympgl}, we can choose
   \[
   u' \in \ker\bigl(\ch{\sympgl} \to \ch{\cycmu}\bigr)\invert
   \]
mapping to $u$ in $\ch{\toruspgl}$. By Proposition~\ref{prop:localization}, we can write
   \[
   u' = v + \chern_{p-1}(D_{p})w,
   \]
where $v$ is in $\ch{\pgl}\invert$ and $w$ is in $\ch{\sympgl}\invert$. The image of $\chern_{p}(D_{p})$ in $\ch{\toruspgl}$ is $0$, because $\toruspgl$ acts trivially on $D_{p}$; so the image of $v$ in $\ch{\toruspgl}$ equals $u$. But there is no reason why $v$ should map to $0$ in $\ch{\cycmu}$.

Let us denote by $\overline{v}$ and $\overline{w}$ the images of $v$ and $w$ respectively in
   \[
   \ch{\cycmu}\invert = \ZZ[1/(p-1)!][\xi,\eta]/(p\xi,p\eta);
   \]
the restriction of $\chern_{p-1}(D_{p})$ equals $-\xi^{p-1}$, so we have $\overline{v} - \xi^{p-1}\overline{w} = 0$. On the other hand $\overline{v}$ is contained in 
   \[
   \bigl(\ch{\cycmu}\bigr)^{\slfinite}\invert = \ZZ[1/(p-1)!][q,r]/(pq,pr);
   \]
since $\overline{v}$ is contained in the ideal of $\ZZ[1/(p-1)!][\xi,\eta]/(p\xi,p\eta)$ generated by $\xi$, and the images of $q$ and $r$ in
   \[
   \ZZ[1/(p-1)!][\xi,\eta]/(\xi,p\eta) = \ZZ[1/(p-1)!][\eta]/(p\eta)
   \]
are $\eta^{p^{2}-p}$ and $0$, we see that $\overline{v}$ is a multiple of $r$; hence we can write $\overline{v} = r\phi(q,r)$, where $\phi$ is a polynomial with coefficients in $\ZZ[1/(p-1)!]$. Set
   \[
   v' = v - \rho\phi(-\chern_{p^{2}-p},\rho);
   \]
then $v'$ restricts to $0$ in $\ch{\cycmu}$, and its image in $\ch{\toruspgl}$ equals the image of $v$, which is $u$, because $\rho$ maps to $0$.

This concludes the proof of Lemma~\ref{lem:ontoker-pgl}.
\end{proof}

Set $K = \ker\bigl(\ch{\pgl} \to \ch{\cycmu}\bigr)$ and $L =  \ker\bigl(\invtorus \to \ch{\mmu}\bigr)$. The induced homomorphism $K \arr L$ is an isomorphism, according to Lemma~\ref{lem:ontoker-pgl}.

Consider the subring $\mathbb{Z} \oplus L \subseteq \invtorus$; Proposition \refall{prop:describe-cycpgl}{5} gives us a copy $\mathbb{Z} \oplus K$ of it inside $\ch{\pgl}$. To finish the proof of Theorem~\ref{thm:main-splitting} we need to extend this splitting to all of $\invtorus$. According to Proposition~\ref{prop:ker-whole}, we have that $\invtorus$ is generated as an algebra over $\mathbb{Z} \oplus L$ by the single element $\delta$. We need to find a lifting for $\delta$; this is provided by the following lemma.

\begin{lemma}
The restriction of $\chern_{p^2-p}(\liesl) \in \ch{\pgl}$ to $\ch[p^2-p]{\toruspgl}$ equals $\delta$.
\end{lemma}

\begin{proof}
We use the notation in the beginning of Section~\ref{sec:mainresults}. The representations $\liesl$ and $\liegl = \liesl \oplus k$ of $\pgl$ have the same Chern classes. If $V = k^n$ is the standard representation of $\GL$, then $\liegl = V \otimes V^\vee$ has total Chern class
   \[
   \chern_t(\liegl) = \prod_{i, j}\bigl(1 + t(x_i-x_j)\bigr)
   = \prod_{i \neq j}\bigl(1 + t(x_i-x_j)\bigr)
   \]
in $\ch{\torusgl}$; but $\ch{\toruspgl} \subseteq \ch{\torusgl}$, so the thesis follows.
\end{proof}

Set $\delta_1 = \chern_{p^2-p}(\liesl) \in \ch{\pgl}$. We consider the subring $(\mathbb{Z}\oplus K)[\delta_1]$ of $\ch{\pgl}$; to finish the proof of the theorem we have left to show that it maps injectively into $(\mathbb{Z}\oplus L)[\delta] = \invtorus$.

Let us take a homogeneous element $x \in (\mathbb{Z}\oplus K)[\delta_1]$ that maps to $0$ in $\invtorus$; according to Proposition~\ref{prop:keyinjectivity}, to check that it is $0$ it is enough to prove that it restricts to $0$ into $\ch{\cycmu}$. Write
   \[
   x = a_0 + a_1 \delta_1 + a_2 \delta_1^2 + a_3 \delta_1^3 + \cdots.
   \]
The $a_i$ of positive degree are in $K$, and therefore map to $0$ in $\ch{\cycmu}$ by definition; so there can be at most one term that does not map to zero, and that has to be of the form $h \delta_1^d$, where $h$ is an integer. However, the restriction of $x$ to $\ch{\mmu} = \mathbb{Z}[\eta]/(p\eta)$ is zero, and since $\delta_1$ restricts to a nonzero multiple of $\eta^{p^2 - p}$ we see that $h$ must be divisible by $p$. This proves that $h \delta_1^d$ also restricts to $0$ in $\ch{\cycmu}$, and completes the proof of the theorem.

\begin{remark}
The splitting $\bigl(\ch{\toruspgl}\bigr)^{\s} \arr \ch{\pgl}$ that we have constructed is not compatible with the splitting $\bigl(\ch{\toruspgl}\bigr)^{\cyc} \arr \ch{\cycpgl}$ constructed in Section~\ref{sec:cycpgl}, in the sense that the diagram
   \[
   \xymatrix{{}\bigl(\ch{\toruspgl}\bigr)^{\s} \ar[r] \ar[d]
       & {}\ch{\pgl} \ar[d] \\
   \bigl(\ch{\toruspgl}\bigr)^{\cyc}\ar[r]
      & \ch{\cycpgl}\hsmash{,}
   }
   \]
where the rows are the splittings and the columns are restrictions, does \emph{not} commute.


\end{remark}

\section{The proofs of the main Theorems}

Let us prove Theorem~\ref{thm:main-ch}. 

First of all, let us check that $\rho$ generates $\ch{\pgl}$ as an algebra over $\invtorus$. Take a homogeneous element $\alpha \in \ch{\pgl}$. The image of $\delta \in \invtorus$ in $\ch{\pgl}$ is $\chern_{p^{2}-p}(\liesl)$, by construction; and this maps to $-q$ in $\ch{\cycmu}$, by Lemma~\ref{lem:restrict-chern}. So there is a polynomial $\phi(x,y)$ with integer coefficients such that $\alpha - \phi(\delta,\rho)$ is in the kernel of the restriction homomorphism $\ch{\pgl} \arr \ch{\cycmu}$; but this kernel is in the image of $\invtorus$, again by construction; and this shows that $\ch{\pgl}$ is generated by $\rho$ as an algebra over $\invtorus$.

The relations given in the statement are satisfied. We have $p\rho = 0$ by construction. Furthermore, by construction the splitting $\invtorus \arr \ch{\pgl}$ sends the kernel of the homomorphism $\invtorus \arr \ch{\mmu}$ into the kernel of $\ch{\pgl}\arr \ch{\cycmu}$; hence, if $u \in \ker(\invtorus \arr \ch{\pgl})$ we have that $u\rho \in \ch{\cycmu}$ goes to $0$ in $\ch{\toruspgl}$, because $\rho$ is torsion, and to $\ch{\cycmu}$. Hence $u\rho = 0$, because of Proposition~\ref{prop:keyinjectivity}.

Let $x$ be an indeterminate, $I$ the ideal in the polynomial algebra 
   \[
   \invtorus{[x]}
   \]
generated by $px$ and by the polynomials $ux$, where $u$ is in the kernel of the restriction homomorphism $\invtorus \arr \ch{\mmu}$; we need to show that the homomorphism $\invtorus{[x]}/I \arr \ch{\pgl}$ that sends $x$ to $\rho$ is an isomorphism. Pick a polynomial $\phi \in \invtorus{[x]}$ such that $\phi(\delta) = 0$ in $\ch{\pgl}$. After modifying it by an element of $I$, we may assume that is of the form $\alpha + \psi(\delta,\rho)$, where $\alpha$ is in the the kernel of $\invtorus \arr \ch{\mmu}$, while $\psi$ is a polynomial in two variables with coefficients in $\FF_{p}$. Since the images of $\delta$ and $\rho$ in $\ch{\cycmu}$, that are $q$ and $r$, are linearly independent in $\FF_{P}[\xi,\eta]$, we see that $\psi$ must be $0$. Hence $\alpha = 0$ in $\ch{\pgl}$; but since $\invtorus$ injects inside $\ch{\pgl}$, we have that $\phi(x) = 0$, as we want.

Next we prove Theorem~\ref{thm:main-H}. We start by proving Corollary~\ref{cor:isom-even}, that says that the cycle homomorphism $\ch{\pgl} \arr \H[\even]{\pgl}$ is an isomorphism.

Call $K$ and $L$, respectively, the kernels of the restriction homomorphisms $\ch{\pgl} \arr \bigl(\ch{\cycmu}\bigr)^{\slfinite}$ and $\H[\even]{\pgl} \arr \bigl(\H[\even]{\cycmu}\bigr)^{\slfinite}$; we have a commutative diagram
   \[
   \xymatrix{
   0 \ar[r] & K \ar[r]\ar[d]
      & \ch{\pgl} \ar[r]\ar[d]
      & \bigl(\ch{\cycmu}\bigr)^{\slfinite}\ar[r]\ar[d] & 0\\
   0 \ar[r] & L \ar[r]
     & \H[\even]{\pgl} \ar[r]
     & \bigl(\H[\even]{\cycmu}\bigr)^{\slfinite}\ar[r] & 0\\
   }
   \]
with exact rows. The right hand column is an isomorphism, because of Propositions \ref{prop:describe-ch-cycmu} and \ref{prop:H-invariants-cycmu}. The group $L$ injects into
   \[
   \ker\Bigl(\bigl(\H{\toruspgl}\bigr)^{\s}\to \H{\mmu}\Bigr) = 
   \ker\Bigl(\bigl(\ch{\toruspgl}\bigr)^{\s}\to \ch{\mmu}\Bigr),
   \]
because of Proposition~\ref{prop:keyinjectivity}; on the other hand the restriction homomorphism
   \[
   K \arr \ker\bigl(\bigl(\ch{\toruspgl}\bigr)^{\s}\to \ch{\mmu}\bigr)
   \]
is an isomorphism, because of Lemma~\ref{lem:ontoker-pgl}. This proves that $K \arr L$ is an isomorphism, and this proves Corollary~\ref{cor:isom-even}.

To show that  $\rho$ and $\beta$ generate $\H{\pgl}$ as an algebra over $\invtorus$, take a homogeneous element $\alpha \in \ch{\pgl}$. The element $\rho$ generates $\H[\even]{\pgl}$, because of Theorem~\ref{thm:main-ch} and the fact above.

If $\alpha$ is a homogeneous element of odd degree in $\H{\pgl}$, its image in $\H[\odd]{\cycmu}$ can be written in the form $\phi(q,r)s$, where $\phi$ is an integral polynomial, by Proposition~\ref{prop:H-invariants-cycmu}. Then $\alpha - \phi(-\delta,\rho)\beta$ maps to $0$ in $\H[\odd]{\cycmu}$. On the other hand $\H[\odd]{\pgl}$ injects into $\H[\odd]{\cycmu}$, by Proposition~\ref{prop:keyinjectivity}, and this completes the proof that $\rho$ and $\beta$ generated.

To prove that the given relations generated the ideal of relations is straightforward, and left to the reader.

Finally, let us prove Theorem~\ref{thm:additive-structure}.

Since the homomorphisms $\ch{\pgl}\otimes \QQ \arr \ch{\SL}\otimes\QQ$ and $\ch{\pgl}\otimes \QQ \arr \ch{\SL}\otimes\QQ$ are isomorphisms, the ranks of $\ch[i]{\pgl}$ and $\H[i]{\pgl}$ equal the ranks of $\ch[i]{\SL}$ and $\H[i]{\SL}$. The ranks of the $\H[i]{\pgl}$ are $0$ when $i$ is odd; while for any $m \geq 0$ the rank of $\ch[m]{\pgl} \simeq \H[2n]{\pgl}$ equal the number of monomials of degree $m$ in $\sigma_{2}$, \dots,~$\sigma_{p}$. Such a monomial $\sigma_{2}^{d_{2}} \dots p^{d_{p}}$ can be identified with a partition $\generate{2^{d_{2}} \dots p^{d_{p}}}$ of $m$, so this rank is the numbero of partitions of $m$ with numbers between $2$ and $p$. 

 On the other hand it follows from Theorem~\ref{thm:main-ch} that the torsion part of $\ch{\pgl}$ is a vector space over the field $\FF_{p}$, with a basis given by the elements $\delta^{i}\rho^{j}$, where $i \geq 0$ and $j > 0$. Similarly, from Corollary~\ref{cor:isom-even} we see that the same elements form a basis for $\H[\even]{\pgl}$, while $\H[\odd]{\pgl}$ is an $\FF_{p}$-vector space with a basis formed by the elements $\delta^{i}\rho^{j}\beta$, where $i \geq 0$ and $j \geq 0$.
 
The theorem follows easily from these facts.

\section{On the ring $\invtorus$}\label{sec:invtorus}

If $T$ is a torus, we denote by $\widehat{T}$ the group of characters $T \arr \gm$. We have a homomorphism of $\widehat{T}$ into the additive group $\ch{T}$ that sends each character into its first Chern class: and this induced an isomorphism of the symmetric algebra $\sym_{\ZZ}\widehat{T}$ with $\ch{T}$.

In this section we study the ring of invariants $\invtorus$. It is convenient to view $\invtorus$ as a subring of $\bigl(\ch{\torusgl}\bigr)^{\s}$; this last ring is generated by the symmetric functions $\sigma_{1}$, \dots,~$\sigma_{p}$ of the first Chern characters $x_{1}$, \dots,~$x_{p}$ of the projections $\torusgl \arr \gm$.

If we tensor $\invtorus$ with $\ZZ[1/p]$, then we get a polynomial ring; and it is easy to exhibit generators. The homomorphism of groups of characters
   \[
   \chartorus{\pgl} \arr \chartorus{\SL}
   \]
induced by the projection $\torussl \arr \toruspgl$ is injective, with cokernel $\ZZ/p\ZZ$; hence it becomes an isomorphism when tensored with $\ZZ[1/p]$. Hence
   \[
   \invtorus\invertp \arr \bigl(\ch{\torussl}\bigr)^{\s}\invertp
   \]
is an isomorphism.

According to Lemma~\ref{lem:permutation-rep}, the ring $\bigl(\ch{\torussl}\bigr)^{\s}$ is a quotient
   \[
   \bigl(\ch{\torusgl}\bigr)^{\s}/(\sigma_{1})
   = \ZZ[\sigma_{1}, \dots, \sigma_{p}]/(\sigma_{1})
   = \ZZ[\tau_{2}, \dots, \tau_{p}],
   \]
where we have denoted by $\tau_{i}$ the image of $\sigma_{i}$ in $\ch{\torussl}$. One way to produce elements of $\invtorus$ is to write down explicitly the elements corresponding to the $\sigma_{i}$ in the isomorphism
   \[
   \invtorus\invertp \simeq \ZZ[1/p][\sigma_{2}, \dots, \sigma_{p}]
   \]
and then clear the denominators. 

The composite
   \[
   \xymatrix@R=10pt{
   \ch{\toruspgl}\invertp \ar@{^(->}[r]
      & \ch{\torusgl}\invertp  \ar[r] \ar@{=}[d]
      & \ch{\torussl}\invertp  \ar@{=}[d] \\
      & \ZZ[1/p][x_{1}, \dots, x_{p}] \ar[r]
      &  \ZZ[1/p][x_{1}, \dots, x_{p}]/(\sigma_{1})
   }
   \]
is an isomorphism, and the inverse $\ZZ[1/p][x_{1}, \dots, x_{p}]/(\sigma_{1}) \arr \ch{\toruspgl}\invertp$ is obtained by sending $x_{i}$ to $x_{i  } - \frac{1}{p}\sigma_{1}$. We need to compute the image of the $\sigma_{k}$ in $\ch{\pgl}\invertp \subseteq \ZZ[1/p][\sigma_{1}, \dots, \sigma_{p}]$, and this is given by the following formula (the one giving the Chern classes of the tensor product of a vector bundle and a line bundle).
\begin{lemma}
If $t$ is an indeterminate, we have
   \begin{align*}
      \sigma_{k}(x_{1}+ t, \dots, x_{p} + t) &=
   \sum_{i=0}^{k} \binom{p-k+i}{i} t^{i}\sigma_{k-i}\\
   &= \sigma_{k} 
      + (p-k+1)t \sigma_{k-1}
      + \binom{p-k+2}{2}t^{2}\sigma_{k-2}\\
      &\qquad +\dots 
      + \binom{p-1}{k-1}t^{k-1}\sigma_{1}
      + \binom{p}{k}t^{k}.
   \end{align*}in $\ZZ[x_{1}, \dots, x_{p}, t]$, for $k=0$, \dots,~$p$.
\end{lemma}

\begin{proof}
This follows by comparing terms of degree $k$ in the equality
   \begin{align*}
   \sum_{i=0}^{p}(1+t)^{i}\sigma_{p-i} &= \prod_{i=1}^{p}(1+t+x_{i})\\
   & = \sum_{i=0}^{p}\sigma_{i}(x_{1}+ t, \dots, x_{p} + t).\qedhere
   \end{align*}
\end{proof}

If we subtitute $-\frac{1}{p}\sigma_{1}$ for $t$ we obtain the images of the $\tau_{k}$ in $\invtorus\invertp$; we denote them by $\gamma'_{k}$. In order to get elements of $\invtorus$, we can clear the denominators in the $\gamma'_{i}$; by a straightforward calculation we can check that
   \begin{align*}
   \gamma_{k} &= p^{k-1}\gamma'_{k} \\
   &= \sum_{i=0}^{k-2}(-1)^{i}p^{k-i-1}
      \binom{p-k+i}{i} \sigma_{k-i}\sigma_{1}^{i}
      +(-1)^{k-1}\frac{k-1}{k}\binom{p-1}{k-1}\sigma_{1}^{k}
   \end{align*}
for $k = 2$, \dots,~$p-1$, while
   \begin{align*}
   \gamma_{p} &= p^{p}\gamma'_{p} \\
   &= \sum_{i=1}^{p-2}(-1)^{i}p^{p-i}\sigma_{p-i}\sigma_{1}^{i}
      +(p-1)\sigma_{1}^{p}.
   \end{align*}

From the discussion above we get that $\invtorus\invertp$ is a polynomial ring over $\ZZ[1/p]$ over $\gamma_{2}$, \dots,~$\gamma_{p}$. However, the $\gamma_{i}$ cannot generate $\invtorus$ integrally, because all of them are in the kernel of the homomorphism $\invtorus \arr \ch{\mmu}$, while $\delta\in \bigl(\ch[p^{2} - p]{\toruspgl}\bigr)^{\s}$ is not.

When $p = 3$ the situation is simple. The following result was proved by Vezzosi.

\begin{theorem}[\mbox{\cite[Lemma~3.2]{vezzosi-pgl3}}]\hfil
   \[
   \invtorus[3] = \ZZ[\gamma_{2}, \gamma_{3}, \delta]/
   \bigl(27\delta - 4\gamma_{2}^{3} - \gamma_{3}^{2}\bigr).
   \]
\end{theorem}

\begin{proof}
What follows is essentially the argument given in the proof of Lemma~3.2 in \cite{vezzosi-pgl3}. We have
   \begin{align*}
   \gamma_{2} &= 3\sigma_{2} - \sigma_{1}^{2}\\
   \intertext{and}
   \gamma_{3} &= 27\sigma_{3} - 9\sigma_{1}\sigma_{2} + 2\sigma_{1}^{3}.
   \end{align*}
Let us express $\delta$ as a rational polynomial in $\gamma_{2}$ and $\gamma_{3}$. This is most easily done after projecting into
   \[
   \ch{\torussl} = \ZZ[x_{1}, x_{2}, x_{3}]/(\sigma_{1})
   \]
since we know that $\ch{\toruspgl}$ injects inside $\ch{\torussl}$. Since $\delta$ is the opposite of the classical discriminant $\prod_{1\leq i< j \leq3}^{3}(x_{i}-x_{j})^{2}$, the image of $\delta$ in $\ZZ[x_{1}, x_{2}, x_{3}]/(\sigma_{1})$ equals $4\sigma_{2}^{3} + 27\sigma_{3}^{2}$. The images of $\gamma_{2}$ and $\gamma_{3}$ in $\ZZ[x_{1}, x_{2}, x_{3}]/(\sigma_{1})$ are $3\sigma_{2}$ and $27\sigma_{3}$; hence we get the formula
   \[
   27\delta = 4\gamma_{2}^{3} + \gamma_{3}^{2}
   \]
showing that the relation in the statement of the theorem holds.

We will show that if $\phi \in \invtorus[3] \subseteq \ZZ[\sigma_{1},\sigma_{2},\sigma_{3}]$ is such that $3\phi$ is in $\ZZ[\gamma_2, \gamma_{3}, \delta]$, then $\phi$ is also in $\ZZ[\gamma_2, \gamma_{3}, \delta]$. This implies that $\ZZ[\gamma_2, \gamma_{3}, \delta] = \invtorus[3]$, because we know that $\gamma_{2}$ and $\gamma_{3}$ generate $\invtorus[3]\otimes \ZZ[1/3]$, hence if $\phi \in \invtorus[3]$ then $3^{n}\phi \in \ZZ[\gamma_2, \gamma_{3}, \delta]$ for sufficiently large $n$, and we can proceed by descending induction on $n$.

Write
   \[
   3\phi = p(\gamma_{2}, \gamma_{3}, \delta)
   \in \ZZ[\sigma_{1}, \sigma_{2},\delta] \subseteq \invtorus[3];
   \]
for an integral polynomial $p$. The image of $p(\gamma_{2}, \gamma_{3}, \delta)$ in the polynomial ring $\FF_{3}[\sigma_{1}, \sigma_{2}, \sigma_{3}]$ is $0$; the images of $\gamma_{2}$, $\gamma_{3}$ and $\delta$ in $\FF_{3}[\sigma_{1}, \sigma_{2}, \sigma_{3}]$ are $-\sigma_{1}^{2}$, $-\sigma_{1}^{3}$ and $\sigma_{2}^{3}$ respectively; and the ideal of relations between these three polynomials is generated by $(-\sigma_{2})^{3} + (-\sigma^{3})^{2}$, which is the image in $\FF_{3}[\sigma_{1}, \sigma_{2}, \sigma_{3}]$ of
   \[
   4\gamma_{2}^{3} + \gamma_{3}^{2} = 27\delta.
   \]
Hence there are two integral polynomials $q$ and $r$  such that we can write
   \begin{align*}
   3\phi &= p(\sigma_{2}, \sigma_{3}, \delta)\\
   &= 3q(\sigma_{2}, \sigma_{3},\delta)
      + 27\delta \cdot r(\sigma_{2}, \sigma_{3},\delta).
   \end{align*}
Dividing by $3$ we see that $\phi$ is in $\ZZ[\sigma_{2}, \sigma_{3}, \delta] \subseteq \invtorus[3]$; and this concludes the proof that $\gamma_{2}$, $\gamma_{3}$ and $\delta$ generate.

Consider the surjective ring homomorphism
   \[
   \ZZ[x_{2}, x_{3}, y] \arr \invtorus[3]
   \]
that sends $x_{i}$ to $\gamma_{i}$ and $y$ to $\delta$; call $I$ its kernel. We know that
   \[
   (27y - 4x_{2}^{3} - x_{3}^{2}) \subseteq I.
   \]
After tensoring with $\ZZ[1/3]$, both rings
   \[
   \ZZ[x_{2}, x_{3}, y]/(27y - 4x_{2}^{3} - x_{3}^{2})
   \quad\text{and}\quad
   \invtorus[3]
   \]
become polynomial rings $\ZZ[1/3][x_{2}, x_{3}]$; hence, if $f \in I$, some multiple $3^{n}f$ is in $(27y - 4x_{2}^{3} - x_{3}^{2})$. But this implies that $f$ is in $(27y - 4x_{2}^{3} - x_{3}^{2})$ because $3$ is a prime in the unique factorization domain $\ZZ[x_{2}, x_{3}, y]$ and does not divide $27y - 4x_{2}^{3} - x_{3}^{2}$, so
   \[
   (27y - 4x_{2}^{3} - x_{3}^{2}) = I,
   \]
as claimed.
\end{proof}

From this, Theorem~\ref{thm:pgl3} follows easily.

As $p$ grows, the calculations become very complicated very quickly. The obvious generalization of the result above, that $\invtorus$ is generated by the $\gamma_{i}$ and $\delta$, fails badly.  When $p$ is larger than $3$, it is not hard to see that they fail to generate already in degree $4$. When $p = 5$ the ring $\invtorus[5]$ has $9$ generators, in degrees $2$, $3$, $4$, $5$, $6$, $7$, $9$, $12$ and $20$; with some pain, it is possible to write them down explicitly. The generators in degree $2$ and $3$ are $\gamma_{2}$ and $\gamma_{3}$. With more work it shoud also be possible to find the relations among them.

There are other approaches to calculations other than the one given here for $p = 3$; but none of them seem to give a lot of information in the general case.

\nocite{edidin-graham-equivariant}
\nocite{gottlieb-euler}
\nocite{grothendieck-brauer3}
\nocite{kono-mimura-shimada}
\nocite{kono-mimura}
\nocite{toda-classifying}

\bibliographystyle{amsalpha}
\bibliography{mrabbrev,VistoliRefs}

\end{document}